\documentclass[10pt,a4paper,final]{article}
\usepackage{url}
\usepackage[a4paper,portrait]{geometry}
\usepackage{amsmath,amsfonts,latexsym,amssymb,amsthm}

\newcommand{\EM}{\ensuremath}

\makeatletter
\newcommand{\@THMSTYLES}{%
  \newtheoremstyle{bodyrm}
  {3pt}
  {3pt}
  {}
  {}
  {\bfseries\sffamily}
  {.}
  { }
  {}
  \newtheoremstyle{bodyit}
  {3pt}
  {3pt}
  {\itshape}
  {}
  {\bfseries\sffamily}
  {.}
  { }
  {}
}
\newcommand{\THMEN}{%
  \@THMSTYLES
  \theoremstyle{bodyit}
  \newtheorem{thm}{Theorem}[section]%
  \newtheorem{cor}[thm]{Corollary}%
  \newtheorem{prop}[thm]{Proposition}%
  \newtheorem{lem}[thm]{Lemma}%
  \theoremstyle{bodyrm}%
  \newtheorem{defi}[thm]{Definition}%
  \newtheorem{xpl}[thm]{Example}%
  \newtheorem{exo}[thm]{Exercise}%
  \newtheorem{hyp}[thm]{Hypothesis}%
  \newtheorem{eur}[thm]{Heuristics}%
  \newtheorem{pro}[thm]{Problem}%
  \newtheorem{rem}[thm]{Remark}%
  \newtheorem{prp}[thm]{Property}%
}
\newcommand{\THMFR}{%
  \@THMSTYLES
  \theoremstyle{bodyit}
  \newtheorem{thm}{Théorème}[section]%
  \newtheorem{cor}[thm]{Corollaire}%
  \newtheorem{prop}[thm]{Proposition}%
  \newtheorem{lem}[thm]{Lemma}%
  \theoremstyle{bodyrm}%
  \newtheorem{defi}[thm]{Définition}%
  \newtheorem{xpl}[thm]{Exemple}%
  \newtheorem{hyp}[thm]{Hypothèse}%
  \newtheorem{rem}[thm]{Remarque}%
}
\makeatother

\makeatletter
\newcommand{\SMALLSECS}{%
 \renewcommand{\section}{\@startsection%
  {section}
  {1}
  {0em}
  {\baselineskip}
  {0.5\baselineskip}
  {\normalfont\large\bfseries}}
 \renewcommand{\subsection}{\@startsection%
  {subsection}
  {2}
  {0em}
  {\baselineskip}
  {0.25\baselineskip}
  {\normalfont\bfseries}}
}
\makeatother

\makeatletter
\providecommand{\timenow}{\@tempcnta\time
\@tempcntb\@tempcnta
\divide\@tempcntb60
\ifnum10>\@tempcntb0\fi\number\@tempcntb
\multiply\@tempcntb60
\advance\@tempcnta-\@tempcntb
:\ifnum10>\@tempcnta0\fi\number\@tempcnta}
\makeatother

\makeatletter
\newcommand{\versiondetravail}{%
 \renewcommand{\@evenfoot}{%
 \hfil{\tiny\texttt{%
   Version préliminaire, compilée le \today{} à \timenow.}\hfill}}%
 \renewcommand{\@oddfoot}{\@evenfoot}%
}
\makeatother

\newcommand{\dB}{\EM{\mathbb{B}}}

\newcommand{\dE}{\EM{\mathbb{E}}}

\newcommand{\dH}{\EM{\mathbb{H}}}

\newcommand{\dN}{\EM{\mathbb{N}}}

\newcommand{\dR}{\EM{\mathbb{R}}}
\newcommand{\dS}{\EM{\mathbb{S}}}

\newcommand{\rD}{\EM{\mathrm{D}}}

\newcommand{\rI}{\EM{\mathrm{I}}}

\newcommand{\rL}{\EM{\mathrm{L}}}

\newcommand{\cA}{\EM{\mathcal{A}}}
\newcommand{\cB}{\EM{\mathcal{B}}}
\newcommand{\cC}{\EM{\mathcal{C}}}
\newcommand{\cD}{\EM{\mathcal{D}}}
\newcommand{\cE}{\EM{\mathcal{E}}}
\newcommand{\cF}{\EM{\mathcal{F}}}

\newcommand{\cH}{\EM{\mathcal{H}}}
\newcommand{\cI}{\EM{\mathcal{I}}}

\newcommand{\cL}{\EM{\mathcal{L}}}
\newcommand{\cM}{\EM{\mathcal{M}}}
\newcommand{\cN}{\EM{\mathcal{N}}}

\newcommand{\cP}{\EM{\mathcal{P}}}

\newcommand{\bE}{\EM{\mathbf{E}}}

\newcommand{\bH}{\EM{\mathbf{H}}}

\newcommand{\bJ}{\EM{\mathbf{J}}}

\newcommand{\bL}{\EM{\mathbf{L}}}

\newcommand{\bP}{\EM{\mathbf{P}}}

\newcommand{\bS}{\EM{\mathbf{S}}}

\newcommand{\bW}{\EM{\mathbf{W}}}

\newcommand{\al}{\alpha}
\newcommand{\be}{\beta}

\newcommand{\de}{\delta}
\newcommand{\ga}{\gamma}

\newcommand{\la}{\lambda}

\newcommand{\Om}{\Omega}
\newcommand{\om}{\omega}

\newcommand{\Si}{\Sigma}
\newcommand{\si}{\sigma}
\newcommand{\Te}{\Theta}
\newcommand{\te}{\theta}

\newcommand{\veps}{\varepsilon}

\newcommand{\p}[4]{{#3}\!\left#1{#4}\right#2} 

\newcommand{\ABS}[1]{\EM{{\left| #1 \right|}}} 
\newcommand{\BRA}[1]{\EM{{\left\{#1\right\}}}} 
\newcommand{\DP}[1]{\EM{{\left<#1\right>}}} 
\newcommand{\NRM}[1]{\EM{{\left\| #1\right\|}}} 
\newcommand{\OSC}[1]{\EM{{\p(){\mathrm{osc}}{#1}}}} 
\newcommand{\PAR}[1]{\EM{{\left(#1\right)}}} 
\newcommand{\pd}{\EM{\partial}} 
\newcommand{\SBRA}[1]{\EM{{\left[#1\right]}}} 
\newcommand{\LIP}[1]{\EM{\|#1\|_{\mathrm{Lip}}}} 

\newcommand{\entf}[1]{\mathbf{Ent}_{#1}}
\newcommand{\ent}[2]{\p(){\entf{#1}}{#2}}

\newcommand{\moyf}[1]{\bE_{#1}}
\newcommand{\moy}[2]{\p(){\moyf{#1}}{#2}}

\newcommand{\covf}[1]{\mathbf{Cov}_{#1}}
\newcommand{\cov}[3]{\p(){\covf{#1}}{#2,#3}}

\newcommand{\varf}[1]{\mathbf{Var}_{#1}}
\newcommand{\var}[2]{\p(){\varf{#1}}{#2}}

\newcommand{\Det}[1]{\mathrm{Det}\,}
\newcommand{\TR}[1]{\p(){\mathrm{Tr}}{#1}}

\newcommand{\GA}{\boldsymbol{\Gamma}}
\newcommand{\GD}{\GA_{\!\!{\mathbf 2}}}
\newcommand{\GI}{\bL}
\newcommand{\GR}{\nabla}

\newcommand{\LA}{\boldsymbol{\Delta}}

\newcommand{\SGf}[1]{{\mathbf P}_{#1}}
\newcommand{\SG}[2]{\p(){\SGf{\!#1}}{#2}}

\newcommand{\Hess}[1]{{\p(){\mathrm{Hess}}{#1}}}

\newcommand{\inter}[1]{{\overset{\circ}{#1}}} 

\renewcommand{\leq}{\leqslant}
\renewcommand{\geq}{\geqslant}
\newcommand{\bs}{\EM{\backslash}}

\newcommand{\ds}[1]{\EM{\displaystyle{#1}}}

\SMALLSECS
\theoremstyle{plain}
\newtheorem{thm}{Theorem}[]
\newtheorem{cor}[thm]{Corollary}
\newtheorem{prop}[thm]{Proposition}

\newtheorem{defi}[thm]{Definition}
\theoremstyle{plain} 
\newtheorem{rem}[thm]{Remark}

\title{Entropies, Convexity, and Functional Inequalities\\
 {\small On Phi-entropies and Phi-Sobolev inequalities}}

\newcommand{\mykeywords}{%
 Spectral gap, %
 Poincar{\'e} inequality, %
 Log-Sobolev inequality, %
 Markov Processes, %
 Diffusion Processes, %
 Kullback-Leibler relative entropy, %
 Shannon Entropy, %
 Boltzmann Entropy, %
 Fisher information, %
 L\'{e}vy Processes, %
 Malliavin Calculus, %
 Analysis on Path Space.}

\newcommand{\mysubjclass}{%
 46-99, 
 60J60, 
 26D10, 
 58D25, 
 39B72, 
 58J65, 
 47D07, 
 94A15, 
 94A17
 .}

\author{Djalil \textsc{Chafa\"{i}}}

\date{{\medskip\footnotesize %
  Preprint -- November 26, 2002.\\
  Revised September 29, 2003.\\
  Revised March 12, 2004.}}

\begin{document}

\newcommand{\HYPOT}[1]{{\footnotesize\textbf{(\cH{#1})}}}
\newcommand{\DBLI}{\EM{\cI^{(2)}}}
\newcommand{\whPhi}{\EM{\widehat\Phi}}
\newcommand{\VARF}[2]{\varf{#1}^{#2}}
\newcommand{\VAR}[3]{\p(){\VARF{#1}{#2}}{#3}}
\newcommand{\ENTF}[2]{\entf{#1}^{#2}}
\newcommand{\ENT}[3]{\p(){\ENTF{#1}{#2}}{#3}}
\newcommand{\HENTF}[2]{\bH{#1}^{#2}}
\newcommand{\HENT}[3]{\p(){\HENTF{#1}{#2}}{#3}}
\newcommand{\SENTF}[2]{\bS{#1}^{#2}}
\newcommand{\SENT}[3]{\p(){\SENTF{#1}{#2}}{#3}}
\newcommand{\FISHF}[2]{\bJ_{#1}^{#2}}
\newcommand{\FISH}[3]{\p(){\FISHF{#1}{#2}}{#3}}

\maketitle

\begin{abstract}
  Our aim is to provide a short and self contained synthesis which generalise
  and unify various related and unrelated works involving what we call
  $\Phi$-Sobolev functional inequalities. Such inequalities related to
  $\Phi$-entropies can be seen in particular as an inclusive interpolation
  between Poincar\'{e} and Gross logarithmic Sobolev inequalities.  In
  addition to the known material, extensions are provided and improvements are
  given for some aspects. Stability by tensor products, convolution, and
  bounded perturbations are addressed.  We show that under simple convexity
  assumptions on $\Phi$, such inequalities hold in a lot of situations,
  including hyper-contractive diffusions, uniformly strictly log-concave
  measures, Wiener measure (paths space of Brownian Motion on Riemannian
  Manifolds) and generic Poisson space (includes paths space of some pure
  jumps L\'{e}vy processes and related infinitely divisible laws). Proofs are
  simple and relies essentially on convexity.  We end up by a short parallel
  inspired by the analogy with Boltzmann-Shannon entropy appearing in
  Kinetic Gases and Information Theories.
\end{abstract}

{\footnotesize\noindent
\textbf{Subject Class. MSC-2000} : \mysubjclass
}

\tableofcontents
\section{Introduction}
\label{se:intro}

Let $\Phi:\cI\to\dR$ be a smooth \emph{convex} function defined on a closed
interval $\cI$ of $\dR$ not necessarily bounded.  Let $\mu$ be a positive
measure on a Borel space $(\Om,\cF)$.  The $\Phi$-entropy functional
$\ENTF{\mu}{\Phi}$ is defined on the set of $\mu$-integrable functions
$f:(\Om,\cF)\to(\cI,\dB(\cI))$ by the following formula:
$$
\ENT{\mu}{\Phi}{f} := \int_{\Om}\!\Phi(f)\,d\mu-\Phi\PAR{\int_\Om\!f\,d\mu}.
$$
Obviously, such a formula makes sense only when $\int_\Om\!f\,d\mu\in\cI$,
which is always the case when $\mu$ is a probability measure.  Unless
otherwise stated, the $\Phi$-entropy in the sequel will be always considered
for probability measures. One has then that
\begin{equation}\label{eq:def-phi-ent}
\ENT{\mu}{\Phi}{f}=\moy{\mu}{\Phi(f)}-\Phi\PAR{\moyf{\mu}{f}}. 
\end{equation}
In addition, depending on $\cI$, the integrability condition on $f$ can be
relaxed and Jensen inequality implies that the $\Phi$-entropy functional takes
its values in $\dR_+\cup\{+\infty\}$. Moreover, it is convex with respect to
its functional argument \emph{at fixed mean}.  As we will see, the global
convexity requires more assumptions on $\Phi$.  If $f$ is $\mu$-a.s. constant,
$\ENT{\mu}{\Phi}{f}$ vanishes, and the converse is true when $\Phi$ is
strictly convex. Sometimes, we will drop the $\mu$ subscript in
$\ENTF{\mu}{\Phi}$. For any random variable $X:(\Om',\cF')\to(\Om,\cF)$, we
will denote 
$$
\ENT{}{\Phi}{f(X)}:=\ENT{\cL(X)}{\Phi}{f}.
$$
The classical variance and entropy can be recovered since we have
\begin{equation}\label{eq:std-var-ent}
  \ENTF{\mu}{x\mapsto x^2}=\varf{\mu}
  \text{\quad and \quad}
  \ENTF{\mu}{x\mapsto x\log x}=\entf{\mu}.
\end{equation}
Notice that the $\Phi$-entropy functional $f\mapsto\ENT{\mu}{\Phi}{f}$ is 
neither homogeneous nor translation invariant in general. 
Nevertheless, since it is non-negative and vanishes when its argument $f$ is a 
constant function, it can be a good candidate as a left hand side of a Sobolev 
like functional inequality where the right hand side is a Dirichlet form.

Actually, the term ``$\Phi$-entropy'' is quite arbitrarily chosen, since we
can speak about ``$\Phi$-variance'' too, but this term is perhaps more adapted
to the quantity 
\begin{equation}\label{eq:def-phi-var}
  \VAR{\mu}{\Phi}{f}:=\moy{\mu}{\Phi\PAR{f-\moyf{\mu}{f}}},
\end{equation}
which is translation invariant and gives the classical variance $\varf{\mu}$
when $\Phi(x)=x^2$, but the entropy $\entf{\mu}$ cannot be recovered.  
One can remember the famous quotation from John Von Neumann about \emph{entropy}, which can be
found for example in \cite[Chap. 10]{MR2002g:46132}.  Notice that similar
$\Phi$-entropies appears with that name in a slightly different forms in a lot
of papers related to Information Theory and Convex Analysis fields, see for
example \cite{MR84e:94009}, \cite{MR84e:62085}, \cite{MR84d:94009},
\cite{MR94b:94010} and \cite{MR87k:94007} and references therein. The
$\Phi$-entropy is related to the so called $(h,\Phi)$-entropies, see for
example \cite{MR98c:62008} and references therein. The $\Phi$-entropy is also
known as $J$-divergence ($J$ stands for Jensen). See also
\cite{MR29:1671,MR48:13470} for the similar notion of $\phi$-divergence.

Let $(X_t)_{t\geq 0}$ be a Markov process on a Polish space $\Om$ equipped
with its Borel $\si$-field, say for example $(\dR^d,\dB(\dR^d))$.  We define
the classical associated Markov semi-group $\PAR{\SGf{t}}_{t\geq 0}$ acting on
$\cC_b(\Om,\dR)$ by:
\begin{equation}\label{eq:def-sg}
  \SG{t}{f}(x) := \moy{\mu}{f(X_t)\,\vert\,X_0=x}.
\end{equation}
Let us assume that there exists an invariant measure $\mu$, i.e. a positive
Borel measure $\mu$ on $\Om$ stable by $\PAR{\SGf{t}}_{t\geq 0}$. When $\mu$
is a \emph{probability}, we get that $\cL(X_0)=\mu$ implies $\cL(X_t)=\mu$ for
all $t\in\dR_+$.  We denote by $\GI:=\pd_{t=0^+}\,\SGf{t}$ the infinitesimal
generator of $\PAR{\SGf{t}}_{t\geq 0}$ with domain $\cD(\GI)\in\rL^2(\mu)$,
and by $\GA$ the associated ``carr\'{e} du champ'' operator defined by:
\begin{equation}\label{eq:def-gamma}
  \GA(g,h):=\frac{1}{2}\,\PAR{\GI(gh)-g\,\GI h-h\,\GI g}.
\end{equation}
By convention, $\GA g := \GA(g,g)$.  For Markov processes, $\GA g$ is
always non-negative. The invariance of a probability measure $\mu$ is
equivalent to
$$
\forall f\in\cD(\GI),\quad \moy{\mu}{\GI f}=0.
$$
We say that a positive measure $\mu$ is symmetric if and only if $\GI$ is
symmetric in $\rL^2(\mu)$, i.e.
$\DP{f\,,\,\GI\,g}_{\rL^2(\mu)}=\DP{g\,,\,\GI\,f}_{\rL^2(\mu)}$ for all $f$
and $g$ in $\cD(\GI)$.  Measure $\mu$ is invariant if it is symmetric but the
converse is false in general.  Symmetric measures lead to an integration by
parts formula
\begin{equation}\label{eq:ipp}
  -\DP{g\,,\,\GI\,h}_{\rL^2(\mu)}%
  =\DP{\GA(g,\,h)}_{\rL^2(\mu)}%
  =\DP{\GA(h,\,g)}_{\rL^2(\mu)}.
\end{equation}
The reader may find an introduction to the analysis of Markov semi-groups in
\cite{MR95m:47075} and \cite{bakry-markov-2002} for example, where the
delicate problem of the existence of an algebra of functions $\cA$ stable by
semi-group and generator is addressed. By invariance of $\mu$ and Jensen
inequality for $\Phi$:
\begin{align*}
\ENT{\mu}{\Phi}{\SGf{t}{f}}
 &= \moy{\mu}{\Phi(\SGf{t}{f})}-\Phi(\moyf{\mu}{f}) \\
 &\leq \moy{\mu}{\SGf{t}{\Phi(f)}}-\Phi(\moyf{\mu}{f}) \\
 &= \ENT{\mu}{\Phi}{f}.
\end{align*}
On the other hand, if the semi-group is $\rL^2$-ergodic,
$\ENT{\mu}{\Phi}{\SGf{t}{f}}$ converges to $0$ when $t$ tends to $+\infty$,
and we get that:
$$
0=\ENT{\mu}{\Phi}{\SGf{\infty}{f}}
\leq 
\ENT{\mu}{\Phi}{\SGf{t}{f}}
\leq
\ENT{\mu}{\Phi}{\SGf{0}{f}}=\ENT{\mu}{\Phi}{f}.
$$
Actually, one can show that any $\Phi$-entropy related to the invariant
measure of a Markov process is non-increasing along the associated Markovian
semi-group:

\begin{prop}[DeBruijn like property for Markov semi-groups]
  \label{pr:debruijn}
  Let $(X_t)_{t\geq 0}$ be a Markov process on a Polish space $\Om$ equipped
  with its Borel $\si$-field. Let $\PAR{\SGf{t}}_{t\geq 0}$ be the associated
  Markov semi-group with infinitesimal generator $\GI$ and ``carr\'{e} du
  champ'' $\GA$.  Assume that $\mu$ is an invariant probability measure.
  Then, for any suitable function $f:\Om\to\cI$ and any $t>0$:
  \begin{equation}\label{eq:debruijn-inv}
    \pd_t\,\ENT{\mu}{\Phi}{\SGf{t}{f}} 
    =\moy{\mu}{\Phi'(\SGf{t}{f})\,\GI\SGf{t}{f}} \leq 0,
  \end{equation}
  When $\mu$ is symmetric, one has the following formulation
  \begin{equation}\label{eq:debruijn-sym}
    \pd_t\,\ENT{\mu}{\Phi}{\SGf{t}{f}} 
    =-\moy{\mu}{\GA(\Phi'(\SGf{t}{f}),\,\SGf{t}{f})}.
  \end{equation}
  Moreover, when $\PAR{X_t}_{t\geq 0}$ is a diffusion process, one has:
  \begin{equation}\label{eq:debruijn-diffusion}
  \pd_t \ENT{\mu}{\Phi}{\SGf{t}{f}} 
  =-\moy{\mu}{\Phi''(\SGf{t}{f})\,\GA\SGf{t}{f}}.
  \end{equation}
\end{prop}

\begin{proof}
  Equality in \eqref{eq:debruijn-inv} follows immediately from definition of
  $\ENTF{\mu}{\Phi}$ and $\GI$. The symmetric case \eqref{eq:debruijn-sym}
  comes by integration by parts \eqref{eq:ipp}.  Let us show that the right
  hand side of \eqref{eq:debruijn-inv} is $\leq 0$.  Jensen inequality for
  convex function $\Phi$ and probability measure $\SGf{t}$ yields
  $\Phi(\SG{t}{g}) \leq \SG{t}{\Phi(g)}$ for any $t>0$, and hence
  $\Phi'(g)\,\GI\,g \leq \GI\Phi(g)$. Thus, by invariance of $\mu$, we get
  that $\moy{\mu}{\Phi'(g)\,\GI\,g} \leq 0$, which gives the result when
  $g=\SGf{t}{f}$.  Finally, the diffusion case \eqref{eq:debruijn-diffusion}
  comes from the fact that we have then the so called ``chain rule formula''
  \begin{equation}
  \GA(\al(g),\,h)=\al'(g)\,\GA(g,\,h).\label{eq:chain-rule}
  \end{equation}
  Recall that the operator $\GI$ is a diffusion operator if and only if for
  any $f_1,\ldots,f_k$ in $\cD(\GI)$ and any smooth function $\al:\dR^k\to\dR$
  such that $\al(f_1,\ldots,f_k)\in\cD(\GI)$:
  \begin{equation}
  \GI(\al(f_1,\ldots,f_k)) =\sum_{i=1}^k
  (\pd_i\,\al)(f_1,\ldots,f_k)\,\GI f_i +\!\!\sum_{1\leq i,j\leq k}
  (\pd_{ij}^2\,\al)(f_1,\ldots,f_k)\,\GA(f_i,f_j).\label{eq:def-diffusion}
  \end{equation}
  We have to mention that the diffusion property makes really sense only in
  continuous space settings and implies roughly that $\GI$ is a second order
  linear partial differential operator without constant part, cf.
  \cite{bakry-markov-2002,MR95m:47075}.  Finally, one can observe that the
  convexity of $\Phi$ is needed only in order to give the sign in
  \eqref{eq:debruijn-inv}: the $\Phi$-entropy $\ENTF{\mu}{\Phi}$ is
  non-increasing along the Markov semi-group when $\Phi$ is convex.
  
  Finally, notice that the term $-\Phi\PAR{\moyf{\mu}{f}}$ in the definition
  of $\ENT{\mu}{\Phi}{f}$ plays no role in the non-increasing property along
  the semi-group since by invariance of $\mu$, one has
  $-\Phi\PAR{\moyf{\mu}{\SG{t}{f}}}=-\Phi\PAR{\moyf{\mu}{f}}$. Therefore, one
  can investigate the non-increasing property along the semi-group for generic
  $(h,\Phi)$-entropies defined by $h\PAR{\moy{\mu}{\Phi(f)}}$.
\end{proof}

Property \eqref{eq:debruijn-inv} tells us that any $\Phi$-entropy related to
the invariant measure on a Markov process is non-increasing along the
Markovian semi-group. Actually, an exponential decrease a $\Phi$-entropy along
the semi-group is equivalent to a functional inequality for $\mu$:

\begin{cor}[Exponential decrease of $\ENTF{}{\Phi}$ along a semi-group]%
  \label{co:exp-decr-sg}
  There is an equivalence between:
  \begin{equation}\label{eq:phi-sob}
    \exists\,c\in\dR_+^*,\quad \forall\,f:\Om\to\cI,\quad %
    \ENT{\mu}{\Phi}{f} \leq -c\,\moy{\mu}{\Phi'(f)\,\GI\,f},
  \end{equation}  
  and  
  \begin{equation}\label{eq:exp-decr-sg-phi-ent}
  \exists\,c\in\dR_+^*,\quad \forall\,t\geq 0,\quad \forall f:\Om\to\cI,\quad %
  \ENT{\mu}{\Phi}{\SGf{t}{f}} \leq e^{-t/c }\ \ENT{\mu}{\Phi}{f}.  
\end{equation}
\end{cor}

\begin{proof}
  Obvious from the DeBruijn like property stated in Proposition
  \ref{pr:debruijn} by taking the derivative in $t$.
\end{proof}

Provided that $\cI=\dR$, one can adapt Corollary \ref{co:exp-decr-sg} to the
``$\Phi$-variance'' given by \eqref{eq:def-phi-var}. For Markov processes, the
functional $\FISHF{\mu}{\Phi}$ defined by:
\begin{equation}\label{eq:def-phi-fish-info}
  \FISH{\mu}{\Phi}{f}:=-\int_{\Om}\!\!\Phi'(f)\,\GI\,f\,d\mu(x)
\end{equation}
can be seen as a generalisation of Fisher information, which corresponds to
the cases $\Phi(x)=x\log x$ or $\Phi(x)=x^2$ (both are equivalent for
diffusions due to the chain rule \eqref{eq:chain-rule}). More generally, when
$\cI=\dR_+$, one can define the Kullback-Leibler $\Phi$-entropy (or relative
$\Phi$-entropy) and the $\Phi$-Fisher functionals for any couple of positive
measures $\mu$ and $\nu$ by:
\begin{equation}\label{eq:def-phi-rel-ent}
  \ENT{}{\Phi}{\nu\,\vert\,\mu} =
  \begin{cases}
    \ds{\int_\Om\!\!\whPhi\PAR{\frac{d\nu}{d\mu}}\,d\mu}
    & \text{if $\nu\ll\mu$} \\
    +\infty & \text{if not},
  \end{cases}
\end{equation}
where $\whPhi(u):=\Phi(u)-\Phi(1)\,u$, and
\begin{equation}\label{eq:def-phi-fish-info-mes}
  \FISH{}{\Phi}{\nu\,\vert\,\mu} =
  \begin{cases}
    \ds{-\int_\Om\!\!\Phi'\PAR{\frac{d\nu}{d\mu}}\,%
      \GI\PAR{\frac{d\nu}{d\mu}}}\ d\mu
    & \text{if $\nu\ll\mu$} \\
    +\infty & \text{if not}.
  \end{cases}
\end{equation}
Observe that $\whPhi$ inherits the convexity property from $\Phi$ and that
$(\whPhi(0),\whPhi(1))=(\Phi(0),0)$ and that $(\whPhi)''=\Phi''$: one can
always ``correct'' a convex function by an affine additive part in such a way
that the new function vanishes at a fixed point. If $\nu$ is a probability
measure with $d\nu:=f\,d\mu$, one has
$$
\ENT{}{\Phi}{\nu\,\vert\,\mu} = \ENT{\mu}{\Phi}{f}
\quad\text{and}\quad
\FISH{}{\Phi}{\nu\,\vert\,\mu} = \FISH{\mu}{\Phi}{f}.
$$
It would be interesting to study the role played by such functionals in
Large Deviation Principles, since the ``normal'' entropy appears as a rate
function, cf. \cite{MR90h:60026}. The Kullback-Leibler relative entropy, which
corresponds to $\Phi(x)=x\log x$ appears for example in Sanov Theorem as a
particular convex conjugate functional on probability measures spaces. One can
hope a sort of $\Phi$-Sanov like Theorem involving the $\Phi$-relative
entropy, and the extremal case of the variance is interesting.

\begin{defi}[$\Phi$-Sobolev inequalities]
  In accordance with Corollary \ref{co:exp-decr-sg}, a probability measure
  $\mu$ associated to a Markov process with generator $\GI$ satisfies to a
  $\Phi$-Sobolev inequality of constant $c\in\dR_+$ on the class
  $\cA\subset\{f:\Om\to\cI\}$ if and only if
  \begin{equation}\label{eq:def-phi-sob}
    \forall\,f\in\cA,\ 
    \ENT{\mu}{\Phi}{f} \leq c\,\moy{\mu}{\cE^\Phi(f)},
  \end{equation} 
  where $f\mapsto\cE^\Phi(f)$ is a non-negative ``energy'' functional
  vanishing when $f$ is constant.  The precise choice of such functional will
  depend on the structure on the involved space $\Om$ and measure $\mu$.  When
  $\cE^\Phi(f)=-\Phi'(f)\,\GI{}f$, the ``traditional'' logarithmic Sobolev and
  Poincar\'{e} inequalities correspond respectively to $\Phi(x)=x\log{}x$ and
  to $\Phi(x)=x^2$.  In abstract settings, a $\Phi$-Sobolev inequality for
  $\mu$ takes the form when $\cI=\dR_+$
  \begin{equation}\label{eq:def-phi-sob-mes}
    \forall\,\nu\ll\mu,\ 
    \ENT{}{\Phi}{\nu\,\vert\,\mu} \leq c\,\FISH{}{\Phi}{\nu\,\vert\,\mu}.
  \end{equation}
  Beware that \eqref{eq:def-phi-sob} gives \eqref{eq:def-phi-sob-mes} but the
  converse seems to be true only when $\Phi(x)=x\log x$ (nevertheless, it is
  always true for densities $f$ with respect to $\mu$).  Here again, one can
  prefer the name ``$\Phi$-Poincar\'{e}'' inequalities, but it seems to be an
  illogical choice since log-Sobolev inequalities are not called
  ``log-Poincar\'{e}''! ``Our'' $\Phi$-Sobolev inequalities are close in form
  to the ``Orlicz-Poincar\'{e}'' inequalities mentioned in \cite[p.
  125-126]{guionnet-zegarlinski-ihp-2002}.  The interested reader may find an
  extensive study of similar but quite different inequalities involving
  $\Phi(x)=\ABS{x}^p$ with $p\geq{}1$ in \cite{bobkov-zegarlinski-02} and
  \cite{barthe-cattiaux-roberto}.
\end{defi}

Notice that the two classical $\Phi$-entropies in \eqref{eq:std-var-ent} share
the same special property: $1/\Phi''$ is affine, i.e.  both convex and concave
in the same time. Hence, if one try to generalise these two cases, one can
assume that $1/\Phi''$ is convex. As we will see, it turns out that it is not
a good choice in view of deriving coercive inequalities like Poincar\'{e} or
$\log$-Sobolev inequalities, for which the concavity of $1/\Phi''$ is needed.
Such a condition can be found for example in the unrelated works
\cite{MR2002b:60025}, \cite{MR2001c:60124} and \cite{MR84e:62085}.
Nevertheless, convexity may gives inverse forms of such inequalities. Here are
some possible \emph{additional} assumptions on $\Phi$:
\begin{enumerate}
\item[\HYPOT{1}] \emph{$(u,v)\mapsto\Phi''(u)\,v^2$ is
    non-negative and convex on $\cI\times\cI$ };
\item[\HYPOT{2}] \emph{$(u,v)\mapsto\Phi(u+v)-\Phi(u)-\Phi'(u)\,v$ is non
    negative and convex on $\DBLI$};
\item[\HYPOT{2'}] \emph{$\Phi''$ is convex, non-negative and non-increasing on
    $\cI$};
\end{enumerate}
where $\DBLI:=\{(u,v)\in\dR^2,\,(u,u+v)\in\cI\times\cI\}$. For convenience, we
denote by $\Psi$ the real valued function defined on $\DBLI$ by
\begin{equation}\label{eq:def-psi}
\Psi(u,v):=\Phi(u+v)-\Phi(u)-\Phi'(u)\,v.
\end{equation}
Notice that \HYPOT{1} is equivalent to the convexity of $\Phi$ and $\Phi''$
and $-1/\Phi''$, see Remark \ref{rm:h1-h2} page \pageref{rm:h1-h2}. On the
other hand, \HYPOT{2'} implies \HYPOT{2}, and \HYPOT{2'} is equivalent to
state that $\Phi$, $\Phi''$ and $-\Phi'$ are convex.  Basic examples for both
\HYPOT{1} and \HYPOT{2} and \HYPOT{2'} are given by
\begin{itemize}
 \item $\Phi(x)=x\log x$ on $\cI=\dR_+$;
 \item $\Phi(x)=x^p$ with $1<p<2$ on $\cI=\dR_+$;
 \item $\Phi(x)=x^2$ on $\cI=\dR$.
\end{itemize}
In such examples, the associated $\Phi$-entropy $\ENTF{\mu}{\Phi}$ is
homogeneous. Hypothesis \HYPOT{1} is suitable in continuous settings, whereas
hypothesis \HYPOT{2} and \HYPOT{2'} are useful in discontinuous ones. The
bivariate convexity of $(u,v)\mapsto\Phi''(u)\,v^2$ under \HYPOT{1} and of
$\Psi$ under \HYPOT{2} is the key property to derive $\Phi$-Sobolev
inequalities in the spirit of \eqref{eq:def-phi-sob}. Such functional
inequalities may be investigated in many situations involving Markov
processes:

\begin{center}
{\small
\begin{tabular}[c]{|c|c|c|}\hline
  \textbf{Time} -- \textbf{Space} 
  & \textbf{Continuous space} 
  & \textbf{Discontinuous space} \\ \hline
  \textbf{Continuous time}
  & Diffusions processes 
  & Poisson space and L\'{e}vy processes \\ \hline
  \textbf{Discrete time}
  & \multicolumn{2}{c}{Discrete Markov chains and martingales} \vline\\ \hline
\end{tabular}}
\end{center}
Of course, L\'{e}vy processes are martingales, but it is not the case for all
diffusions. The discrete time case is not really addressed in this work, but
one can consider Random Walks or Bernoulli processes, and some answers can be
found for example in \cite{MR2004c:60043} and references therein. For each
case, one may be interested in inequalities on the whole paths space or more
simply in the law at fixed time (infinite for invariant measure by
ergodicity). There is a lack of chain rule formula \eqref{eq:chain-rule} for
$\GA$ operator in discrete space settings. However, i.i.d. increments can then
help to recover a Brownian Motion like behaviour. In some cases, inequalities
at fixed time can be tensorised to give multi-time inequalities which appears
as inequalities for marginals of the paths space measure. A limiting procedure
can be then used to recover inequalities on paths space (it is the so called
cylindrical method). Roughly, we use here two types of proofs in order to
extract $\Phi$-Sobolev inequalities from ``Markovianity''. The first kind
makes use of the semi-group and the second kind, more powerful, is based on a
martingale representation and gives directly results on paths spaces. Finally,
it is probably possible to establish $\Phi$-Sobolev inequalities under
\HYPOT{1} on some loop spaces, as what is done in \cite{MR99h:58202} and
\cite{MR99f:58222} for the logarithmic Sobolev inequality.

The DeBruijn like property stated in proposition \ref{pr:debruijn} can be seen
as a sort of special case of G{\^a}teau directional derivative of
$\Phi$-entropy like functionals, see for example \cite{MR84e:62085}, where the
$\Phi$-entropy is called $J$-Divergence. Actually, one can state the following
Proposition:

\begin{prop}\label{pr:convex-phi-ent}
  Assume that $\mu$ is a probability measure and that $\Phi$ fulfils
  \HYPOT{1}.  Let $\rL^{1,\Phi}(\mu)$ be the set of measurable functions
  $f:(\Om,\cF)\to\cI$ such that $f\in\rL^1(\mu)$ and $\Phi(f)\in\rL^1(\mu)$.
  Then $\rL^{1,\Phi}(\mu)$ is convex, the $\Phi$-entropy functional
  $f\in\rL^{1,\Phi}(\mu)\mapsto\ENT{\mu}{\Phi}{f}$ is convex, and for any
  $f\in\rL^{1,\Phi}(\mu)$, one has the following duality formula:
  \begin{equation}\label{eq:phi-ent-varfor}
    \ENT{\mu}{\Phi}{f} %
    =\sup_{h\in\rL^{1,\Phi}(\mu)} %
    \BRA{\moy{\mu}{(\Phi'(h)-\Phi'\PAR{\moyf{\mu}{h}})\,(f-h)} %
      +\ENT{\mu}{\Phi}{h}}.
  \end{equation}
\end{prop}

\begin{proof}
  For any couple $f,g$ in $\rL^{1,\Phi}(\mu)$ and any $t\in[0,1]$, one has
  $$
  h_t:=tf+(1-t)g\in\rL^1(\mu),
  $$
  and $\Phi(h_t)\in\rL^1(\mu)$ by convexity of $\Phi$. Thus
  $h_t\in\rL^{1,\Phi}(\mu)$ and $\rL^{1,\Phi}(\mu)$ is convex.  Consider now
  the real function $\al:[0,1]\to\dR_+$ defined by
  $\al(t):=\ENT{\mu}{\Phi}{h_t}$ and let us show that it is convex. Assume
  first that $f$ and $g$ are bounded. An easy computation gives
  $$
  \al'(t) %
  = \moy{\mu}{\Phi'(h_t)(f-g)}-\Phi'(\moy{\mu}{h_t})(\moy{\mu}{f-g}),
  $$
  and
  $$
  \al''(t) %
  = \moy{\mu}{\Phi''(h_t)(f-g)^2} - \Phi''(\moy{\mu}{h_t})(\moy{\mu}{f-g})^2,
  $$
  which is non-negative by virtue of Jensen inequality for the bivariate
  function $(u,v)\mapsto\Phi''(u)\,v^2$. Alternatively, one can use the
  concavity of $1/\Phi''$ and Cauchy-Schwarz inequality.  Therefore, $\al$
  is continuous on $[0,1]$ and convex on $(0,1)$. In other words, for any
  $s,t,\la\in[0,1]$, 
  $$
  \al(\la s + (1-\la)t) \leq \la\al(s)+(1-\la)\al(t).
  $$
  At this stage, the dominated convergence Theorem allows us to drop the
  boundedness assumption on $f$ and $g$. In particular, $\al(\la) \leq
  \la\al(1)+(1-\la)\al(0)$ can be written
  $$
  \ENT{\mu}{\Phi}{\la f+(1-\la)g} %
  \leq \la\ENT{\mu}{\Phi}{f} + (1-\la)\ENT{\mu}{\Phi}{g},
  $$
  and the convexity of $f\in\rL^{1,\Phi}(\mu)\mapsto\ENT{\mu}{\Phi}{f}$ is
  established.  Notice that the convexity of $\Phi$ implies that the
  expression of $\al'$ is valid for any $f$ and $g$ in $\rL^{1,\Phi}(\mu)$.
  Since every convex function is the envelope of its tangents, cf.
  \cite{MR98m:49001,MR2002i:90002}, one gets the following variational
  formula:
  $$
  \ENT{\mu}{\Phi}{f}=\al(1)=\sup_{t\in[0,1]}\BRA{\al(t)+\al'(t)(1-t)},
  $$
  which can be rewritten as \eqref{eq:phi-ent-varfor} since for any fixed $f$ in 
  $\rL^{1,\Phi}(\mu)$,
  $$
  \rL^{1,\Phi}(\mu)=\BRA{tf+(1-t)g,\text{ where } %
    g\in\rL^{1,\Phi}(\mu) \text{ and } t\in[0,1]}.
  $$
  Notice that the value of the sup in \eqref{eq:phi-ent-varfor} is achieved
  for $h=f$.
\end{proof}
The convexity of $\Phi$-entropy like functionals is well known, see for
example \cite{MR84d:94009}, \cite[Thm. 2]{MR84e:62085} and \cite[Lem.
4]{MR2002b:60025}. One can show that this convexity is in some sense
equivalent to \HYPOT{1} via well chosen $f$ and $h$ functions in
\eqref{eq:phi-ent-varfor}. Formally, the Fr\'{e}chet derivatives of this
functional are given by
  \begin{equation*}\label{eq:diff-phi-ent}
    \PAR{D \ENTF{\mu}{\Phi}}(f)(h)
    =\moy{\mu}{\SBRA{\Phi'(f)-\Phi'\PAR{\moyf{\mu}{f}}}\,h},
  \end{equation*}
  and
  \begin{equation*}\label{eq:diff2-phi-ent}
    \PAR{D^2 \ENTF{\mu}{\Phi}}(f)(h,h)
    =\moy{\mu}{\Phi''(f)\,h^2}%
    -\Phi''\PAR{\moyf{\mu}{f}}\,\PAR{\moyf{\mu}{h}}^2.
  \end{equation*}
 We recover the well known formulas for variance and entropy by considering the
appropriate $\Phi$:
\begin{equation}
  \var{\mu}{f} =
  \sup_{h}\BRA{2\,\cov{\mu}{f}{h}-\var{\mu}{h}},\label{eq:var-varfor}
\end{equation}
and
\begin{equation}
  \ent{\mu}{f} = 
  \sup_{h}\BRA{\moy{\mu}{f\log{h}} -\moy{\mu}{f}\log\moy{\mu}{h}}.\label{eq:ent-varfor}
\end{equation}
The duality formula \eqref{eq:phi-ent-varfor} for $\Phi$-entropies can be
found in Convex Analysis and Information Theory literature, at least for
discrete probability measures, see for example \cite{MR87k:94007}. It is quite
straightforward via the convexity of $\Phi$-entropies under \HYPOT{1}. We have
to mention that Pascal Massart \& al gave recently an elementary proof of this
variational formula, cf. \cite{massart-moment-et-al,massart-st-flour}.  See
\cite{MR84e:62085} and \cite{MR98c:62008} for further developments and a
Bayesian and Riemannian Geometry point of view.

\begin{rem}[Bi-convexity of relative $\Phi$-entropy]
  Assume that $\cI=\dR_+$ and let $(\Om,\cF)$ be a measurable space equipped
  by a positive Borel measure $\mu$.  Let $\al$ and $\be$ be two probability
  measures on $(\Om,\cF)$, absolutely continuous with respect to $\mu$, with
  densities $p$ and $q$ respectively.  Then, the relative $\Phi$-entropy
  $\ENT{}{\Phi}{\be\,\vert\,\al}$ defined in \eqref{eq:def-phi-rel-ent} is
  given by
  $$
  \ENT{\mu}{\Phi}{q;p}:=\int_{\Om}\!\whPhi\PAR{\frac{q}{p}}\,p\,d\mu.
  $$
  An easy computation shows that the convexity of $\Phi$ induces the
  bi-convexity of the bivariate function
  $(u,v)\in\dR_+^2\mapsto{}\whPhi(v/u)\,u$. Thus
  $(p,q)\mapsto{}\ENT{\mu}{\Phi}{q;p}$ is also bivariate convex.  Notice that
  at fixed $f$, the $\Phi$-entropy itself $\mu\mapsto{}\ENT{\mu}{\Phi}{f}$ is
  concave!
\end{rem}

\begin{rem}[Convexity and nullity for constants]\label{rm:phimoyf-term}
  One can be surprised by the condition \HYPOT{1} requested for the convexity
  of the $\Phi$-entropy functional $f\mapsto\ENT{\mu}{\Phi}{f}$. Actually, it
  is due to the presence of the $-\Phi(\moyf{\mu}{f})$ term in the definition
  of $\ENTF{\mu}{\Phi}$. If one removes this term, the resulting functional
  $f\mapsto\moy{\mu}{\Phi(f)}$ is always convex since $\Phi$ is convex and
  \HYPOT{1} is not needed.  Nevertheless, in that case, the functional does
  not vanishes when its argument $f$ is constant, and thus it becomes useless
  as the left hand side of a Sobolev like functional inequality.
\end{rem}

\section{Phi-Sobolev for diffusions and log-concave measures}
\label{se:phi-sob-diff-logconc}

Proposition \ref{pr:debruijn} translates the fact that the derivative in time
of the $\Phi$-entropy along a Markovian semi-group can be expressed in terms
of the derivatives in space $\GA$, i.e. Fisher information. As we will see,
for diffusions, a commutation formula between the semi-group and the
derivatives in space yields finally local $\Phi$-Sobolev inequalities. We
learnt the following Theorem for hypercontractive diffusions from \cite[Sect.
3]{MR2001c:60124}, where it is stated in a slightly different manner, see also
\cite{bakry-markov-2002}.

\begin{thm}[$\Phi$-Sobolev inequality for diffusions]
  \label{th:phi-sob-diffusions}
  Let $(\cM,g)$ be a connected complete Riemannian manifold and let
  $\PAR{X_t}_{t\geq 0}$ be a diffusion process on $\cM$ with symmetric
  invariant positive measure $\mu$ absolutely continuous with respect to the
  Riemannian volume measure. If a $\mathrm{CD}(\rho,\infty)$ criterion is
  satisfied with $\rho\geq 0$ (cf. \eqref{eq:curv-dim-crit}, section
  \ref{ss:geom-logconc}, page \pageref{ss:geom-logconc}), then under
  \HYPOT{1}, one has:
  \begin{enumerate}
  \item For any smooth function $f:\cM\to\cI$ and any $t\in\dR_+$:
    \begin{equation}\label{eq:phi-sob-loc-diffusion}
    \ENT{\SGf{t}}{\Phi}{f}\leq \frac{1-e^{-2\rho
        t}}{2\rho}\,\SG{t}{\Phi''(f)\,\GA f},
    \end{equation}
    where $\PAR{\SGf{t}}_{t\geq 0}$ and $\GA$ are as in \eqref{eq:def-sg} and
    \eqref{eq:def-gamma} and where the constant is $t/2$ when $\rho=0$.
  \item If $\rho>0$ and $\mu$ is a probability measure and $(X_t)_{t\geq 0}$
    is $L^2$-ergodic then, for any smooth function $f:\cM\to\cI$:
    \begin{equation}\label{eq:phi-sob-logconc}
      \ENT{\mu}{\Phi}{f} \leq\frac{1}{2\rho}\,\moy{\mu}{\Phi''(f)\,\GA f}.
    \end{equation}
  \end{enumerate}
\end{thm}

\begin{proof}
  The second part may be deduced from the first one via ergodicity by letting
  $t$ tends to $+\infty$. Let us give a direct proof.  Let $\GI$ be the
  infinitesimal generator of the diffusion semi-group.  By the ergodic
  property and Fubini Theorem, it follows that
  \begin{align*}
    \ENT{\mu}{\Phi}{f}
    :&=\moy{\mu}{\Phi(f)}-\Phi(\moyf{\mu}{f})\\
    &=\moy{\mu}{\Phi(\SGf{0}{f})}-\moy{\mu}{\Phi(\SGf{\infty}{f})}\\
    &=-\int_0^\infty\!\! dt\ \moy{\mu}{\pd_t\, \Phi(\SGf{t}{f})}\\
    &=-\int_0^\infty\!\! dt\ \moy{\mu}{\Phi'(\SGf{t}{f})\,\GI\SGf{t}{f}}\\
    &=+\int_0^\infty\!\! dt\ \moy{\mu}{\Phi''(\SGf{t}{f})\,\GA\SGf{t}{f}},
 \end{align*}
 where the last equality is obtained by integration by parts \eqref{eq:ipp}
 which is a consequence of the symmetry of $\mu$ for $\GI$. Now, by the
 diffusion property, the $\mathrm{CD}(\rho,\infty)$ criterion is equivalent to
 the following commutation formula (cf. Section \ref{ss:geom-logconc} page
 \pageref{ss:geom-logconc}):
 \begin{equation}\label{eq:strong-commut-sg}
   \sqrt{\GA\SGf{t}{f}} \leq e^{-\rho t}\ \SGf{t}{\sqrt{\GA f}}.
 \end{equation}
 Therefore, 
 $$
 \Phi''(\SGf{t}{f})\,\GA\SGf{t}{f}
 \leq
 \exp(-2\rho t)\,\Phi''(\SGf{t}{f})\,\PAR{\SGf{t}{\sqrt{\GA f}}}^2,
 $$
 and then by Jensen inequality with the bivariate function $\Phi''(u)\,v^2$ 
 which is convex under \HYPOT{1}:
 $$
 \Phi''(\SG{t}{f})\,\SG{t}{\sqrt{\GA f}}^2
 \leq
 \SG{t}{\Phi''(f)\,\GA f}.
 $$
 Alternatively, one can use the concavity of $1/\Phi''$ and the
 Cauchy-Schwarz inequality.  The desired results follows immediately from
 the invariance of $\mu$: $\moy{\mu}{\SGf{t}{f}}=\moy{\mu}{f}$.  The proof of
 the first part \eqref{eq:phi-sob-loc-diffusion} is quite similar. One just
 has to replace $\mu$ by $\SGf{s}$, and $\SGf{t}$ by $\SGf{t-s}$. The
 integration by parts and the chain rule leading to the expression in $\Phi''$
 must be replaced by the following, relying on the diffusion property of
 $\GI$:
 \begin{align*}
   \pd_s\,\SG{s}{\Phi(\SG{t-s}{f})}
   :&=\SG{s}{\GI \Phi(\SG{t-s}{f})-\Phi'(\SG{t-s}{f})\,\GI\SG{t-s}{f}} \\
   &=\SG{s}{\Phi''(\SG{t-s}{f})\,\GA\SG{t-s}{f}},
 \end{align*}
 where $0<s\leq t$.  Notice that the diffusion property is used directly here,
 which was not the case for the non local inequality, for which it was used
 via the strong commutation property given by the $\mathrm{CD}(\rho,\infty)$
 criterion.  As we will see, there is a lack of the diffusion property in
 discontinuous space settings, and one has to replace
 $(u,v)\mapsto\Phi''(u)\,v^2$ by $\Psi$ which is convex under \HYPOT{2}.
\end{proof}

See \cite{MR2001c:60124} for Gross like hypercontractivity, F.K.G.
inequalities and related aspects of $\Phi$-entropies for diffusion
semi-groups. We have to mention that the diffusion property leads to a more
general result which contains F.K.G. and $\Phi$-Sobolev inequalities as
sub-cases, as explained in \cite[Thm. 4.2]{MR2001c:60124}. Namely, it states
that if $\Phi_1,\Phi_2:\cI\subset\dR\to\dR$ are two smooth functions such that
a particular $5\times 5$ matrix involving the derivatives of $\Phi_1$ and
$\Phi_2$ is positive definite, then for any smooth $f$ and $g$:
\begin{align*}
  \ENT{\mu}{\Phi_1,\,\Phi_2}{f,\,g}
  \leq\ \ &\moy{\mu}{\Phi_1''(f)\,\Phi_2(g)\,\GA f}\\
   &+\moy{\mu}{\Phi_1(f)\,\Phi_2''(g)\,\GA g}\\
   &+2\,\moy{\mu}{\Phi_1'(f)\,\Phi_2'(g)\,\GA(f,g)},
\end{align*}
where the quantity
$$
\ENT{\mu}{\Phi_1,\,\Phi_2}{f,\,g}
:=\moy{\mu}{\Phi_1(f)\,\Phi_2(g)}-\Phi_1(\moy{\mu}{f})\,\Phi_2(\moy{\mu}{g})
$$
is a sort of $(\Phi_1,\,\Phi_2)$-covariance. One has in particular
\begin{equation*}
  \ENTF{\mu}{x\mapsto x,\,x\mapsto x}=\covf{\mu} \quad\text{and}\quad 
  \ENTF{\mu}{x\mapsto 1,\,x\mapsto \Phi(x)}=\ENTF{\mu}{\Phi}.
\end{equation*}

\begin{rem}[Non-negativity of functions]\label{rm:non-neg-func}
  It is possible in some cases to reduce the analysis to non-negative
  functions. Namely, since $\SG{t}{|f|} \geq \ABS{\SG{t}{f}}$, one can write:
  \begin{align*}
    \moy{\mu}{f \GI f} %
    &=\lim_{t\to0^+} \frac{1}{2t}\,\moy{\mu}{f\SG{t}{f}-f^2}\\
    &\leq \lim_{t\to0^+} \frac{1}{2t}\,\moy{\mu}{|f|\SG{t}{|f|}-|f|^2}\\
    &=\moy{\mu}{|f| \GI |f|},
  \end{align*}
  which gives $\moy{\mu}{\GA\ABS{f}} \leq \moy{\mu}{\GA f}$. More generally,
  one can show by the same way that $\moy{\mu}{\Phi''(\ABS{f})\,\GA\ABS{f}}
  \leq \moy{\mu}{\Phi''(\ABS{f})\,\GA f}$. Moreover, one can assume that
  $f\geq\veps>0$ by Fatou Lemma.
\end{rem}

\subsection{Geometry and log-concavity}
\label{ss:geom-logconc}

As we have seen in Theorem \ref{th:phi-sob-diffusions}, a natural framework
for $\Phi$-Sobolev inequalities related to diffusions is Riemannian manifolds.
Let $\GI$ be a Markov generator with symmetric positive measure
$d\mu=\exp(U)\,dv_g$ on a complete connected Riemannian manifold $(\cM,g)$
equipped with its volume measure $v_g$. Function $U:\cM\to\dR$ is taken
smooth.  A basic example is given by the Laplace-Beltrami operator
$\GI=\LA_g+\GR U$ with vector field $\GR U$. Back to the generic case of a
Markov generator $\GI$ on $(\cM,g)$, let us define the iterated functional
quadratic forms $\GA f:=\GA(f,f)$ as is \eqref{eq:def-gamma} and $\GD f :=
\GD(f,f)$ by
\begin{equation}\label{eq:def-gamma2}
  \GD(f,g):=\frac{1}{2}\,\PAR{\GA(f,g)-\GA(f,\GI g)-\GA(g,\GI f)}.
\end{equation}
Then, for any $(\rho,n)\in\dR\times\dN$, the $\mathrm{CD}(\rho,n)$ (or
Bakry-Emery $\GD$) criterion can be expressed as:
\begin{equation}\label{eq:curv-dim-crit}
  \GD\geq\rho\,\GA+\frac{1}{n}\,(\GI)^2,
\end{equation}
where $(\GI)^2 f := (\GI f)^2$. 
For diffusions generators $\GI$, $\mathrm{CD}(\rho,\infty)$ is equivalent to
the strong commutation formula \eqref{eq:strong-commut-sg} between
$\sqrt{\GA}$ and $\SGf{t}$, which is a direct consequence a Bochner
formula\footnote{Actually Bochner-Lichnerowicz-Weitzenboch.}:
$$
\GD(f,f) = (\mathrm{Ric_g}-\GR^2 U)(\GR f,\GR f)+\NRM{\mathrm{Hess} f}_2^2,
$$
see for example \cite{MR99m:60110,MR95m:47075}, \cite{MR2002a:58045} and
\cite[Chap. 5]{MR2002g:46132}. Notice that with $n=\infty$,
\eqref{eq:curv-dim-crit} appears as the infinitesimal form of
\eqref{eq:strong-commut-sg}, via the fact that:
$$
\SGf{t}{f} = f + t\,\GI f + \frac{1}{2}\,t^2\,\GI\GI f + o(t^2),
$$
cf. \cite[Sect. 5.4]{MR2002g:46132}. In absence of the diffusion property,
the Markov semi-group verifies the following weaker commutation formula under
$\mathrm{CD}(\rho,\infty)$:
\begin{equation}\label{eq:weak-commut-sg}
  \GA \SGf{t}{f} \leq e^{-2\rho t}\ \SGf{t}{\GA f}.
\end{equation}
This commutation formula leads to a $\Phi$-Sobolev inequality when
$(u,v)\mapsto\Phi''(u)\,v$ is convex on $\cI\times\dR_+$, which is the case if
and only if $\Phi'''=0$, i.e. $\Phi(u)=au^2+bu+c$ with
$(a,b,c)\in\dR_+\times\dR^2$ (i.e. Poincar\'{e} inequality).

Let us consider the ``simple'' example where $(\cM,g)$ is the standard flat
Euclidean space $(\dR^d,\rI_d)$ and where $\GI=\LA-\GR W \cdot \GR$.  The
symmetric invariant measure $\mu$ is then given by $d\mu(x)=\exp(-W(x))\,dx$
and can be obviously normalised as a probability measure if and only if
$\exp(-W)\in\rL^1(\dR^d,\dR,dx)$. One has $\GA f=\ABS{\GR f}^2$ and
$\mathrm{Ric_g}\equiv 0$ and $\GD f=\NRM{\GR^2 f}_2^2+\GR f \cdot \GR^2 W \GR
f$. Now, if there exists a real number $\rho\geq 0$ such that
$$
\forall x\in\dR^d,\ \GR^2 W(x) \geq \rho \rI_d
$$
as quadratic forms on $\dR^d$, then $\GI$ satisfies
$\mathrm{CD}(\rho,\infty)$. Moreover, if $\rho>0$, measure $\mu$ is uniformly
strictly $\log$-concave and can be normalised as a probability measure.
Obviously, measure $\mu$ can be finite without being uniformly strictly
$\log$-concave. 

\begin{cor}[$\Phi$-Sobolev for uniformly strictly log-concave densities]
  \label{co:phi-sob-logconc-rd}
  Let $\mu$ be a probability measure on $\dR^d$ absolutely continuous with
  respect to Lebesgue measure $dx$ with smooth density. Assume that $\mu$ is
  uniformly strictly log-concave, i.e. that there exists $\rho>0$ such that
  $d\mu(x)=\exp(-H(x))\,dx$ where $H:\dR^d\to\dR$ is smooth and
  $\DP{\Hess{H}(x)\,y;y} \geq \rho\,\NRM{y}_2^2$ for any $x$ and $y$ in
  $\dR^d$. Then under \HYPOT{1} and for any smooth function $f:\dR^d\to\cI$:
  \begin{equation}\label{eq:phi-sob-logconc-Rn}
    \ENT{\mu}{\Phi}{f}%
    \leq\frac{1}{2\,\rho}\,\moy{\mu}{\Phi''(f)\,\ABS{\GR{}f}^2}.
  \end{equation}
  In particular, when $\mu=\cN(m,\Si)$ is a Gaussian measure of mean vector
  $m\in\dR^d$ and positive definite covariance matrix
  $\Si\in\dS\mathrm{ym}_d^{+*}(\dR)$, one has 
  $$
  \rho^{-1}=\max(\si(\Si))
  $$
  where $\si(\Si)$ is the spectrum of $\Si$.  Moreover, the following
  Brascamp-Lieb type inequality holds:
  \begin{equation}\label{eq:gaussian-brascamp-lieb}
    \ENT{\cN(m,\Si)}{\Phi}{f}%
    \leq\frac{1}{2}\,\moy{\cN(m,\Si)}{\Phi''(f)\,\DP{\Si\,\GR{}f;\GR{}f}},
  \end{equation}
  and inequality \eqref{eq:gaussian-brascamp-lieb} remains valid when $\Si$ is
  singular.
\end{cor}
Notice that \eqref{eq:phi-sob-logconc-Rn} is obtained by the comparison
$\Hess{H}\geq \rho\,\rI_d$ to the standard quadratic form
$\DP{\cdot\,;\,\cdot}=\NRM{\cdot}_2^2$ corresponding to the diagonal
covariance matrix $\rI_d$. Such a comparison is essentially unidimensional and
imposes a geometric information loss. Recall that the Brascamp-Lieb inequality
for the strictly log-concave probability measure $\mu$ on $\dR^d$ reads
$$
\forall f:\dR^d\to\dR\text{ smooth},\, %
\var{\mu}{f}\leq \moy{\mu}{\DP{(\Hess{H})^{-1}\GR f;\GR f}},
$$
where $d\mu(x)=\exp(-H(x))\,dx$. When $\mu=\cN(m,\Si)$ with non singular
$\Si$, one has $\Hess{H}=\Si^{-1}$. Unfortunately, to our knowledge, and
except for Gaussian measures, Brascamp-Lieb forms of $\Phi$-Sobolev
inequalities are available only for Poincar\'e inequalities, i.e. when
$\Phi(u)=u^2$, cf. \cite{bobkov-ledoux-bm} for a counter example in the
logarithmic Sobolev case. However, if $\mu$ satisfies to a $\Phi$-Sobolev
inequality on $\dR^d$ with constant $c$ and right hand side
$\moy{\mu}{\Phi''(f)\,\ABS{\GR f}^2}$, one can ask about the best deterministic
symmetric matrix $S$ such that
$$
\ENT{\mu}{\Phi}{f}\leq c\,\moy{\mu}{\Phi''(f)\,\DP{S\GR f;\GR f}}.
$$
Back to the case where $\GI=\LA-\GR\,W\cdot\GR$ on $\dR^d$, an interesting
``degenerated'' case is given by $W(x)=\NRM{x}_r^r$ with $r\in
\{0\}\cup{}(1,2)\cup(2,+\infty)$ since $\GR^2 W(x)$ is then singular at $x=0$.
One has then only $\mathrm{CD}(0,\infty)$, inducing a local version only:
\begin{equation}\label{eq:phi-sob-loc-linear}
  \ENT{\SGf{t}}{\Phi}{f} \leq t\,\SG{t}{\Phi''(f)\,\ABS{\GR f}^2},
\end{equation}
which is the Brownian Motion or (heat semi-group) behaviour corresponding to
the case $r=0$ and $d\mu(x)=dx$. One can expect a better inequality in term of
$t$-dependence of the obtained constant, but the semi-group method seems to
fail since it relies on a commutation formula \eqref{eq:strong-commut-sg}
which is poor when $r\not\in\{0,2\}$.  A perturbative approach (cf.
Proposition~\ref{th:phi-sob-pert} page \pageref{th:phi-sob-pert}) may be used
for the inequalities related to $\mu$ when $r>2$, which corresponds to
$t=+\infty$, but the constant is not sharp in general (cf. \cite[Chap.
6]{MR2002g:46132}).

When $\GR^2 W$ is constant, $(X_t)_{t\geq0}$ is a Brownian Motion or an
Ornstein-Uhlenbeck process, and the commutation formula
\eqref{eq:strong-commut-sg} is exact, i.e. an equality. One can try to
investigate the commutation between $f\mapsto\Phi''(f)\,\GA f$ and $\SGf{t}$,
with a cost related to a sort of convex conjugate $\Phi^*$ of $\Phi$. The idea
is to relate $\Phi$ with a criterion involving $\GA$ and $\GD$ (and thus $W$
when $\GI=\LA-\GR W\cdot\GR$).

\subsection{Interpolation between Poincar\'{e} and logarithmic Sobolev}
\label{ss:interpol-ip-isl}

For simplicity and notational convenience, we restrict ourselves here to
$\dR^d$. Assume that $\Phi:\cI\to\dR$ is convex and smooth and that the
probability measure $\mu$ on $\dR^d$ satisfies to a $\Phi$-Sobolev inequality
of the form:
$$
\exists c>0,\quad\forall f\in\cC_b^\infty(\dR^d,\cI),\quad%
\ENT{\mu}{\Phi}{f} \leq c\,\moy{\mu}{\Phi''(f)\,\ABS{\GR f}^2},
$$
Now, for any $g\in\cC_b^\infty(\dR^d,\cI)$, let us consider
$a\in\inter{\cI}$ such that $\Phi''(a)>0$ and $f_\veps:=a+\veps\,g$ with
$\veps>0$. Then, $f_\veps\in\cC_b^\infty(\dR^d,\cI)$ for $\veps$ sufficiently
small and by Taylor formula for $\Phi$ and the $\Phi$-Sobolev inequality above
for $f_\veps$, one gets when $\veps$ goes to $0^+$:
$$
\var{\mu}{g} \leq 2\,c\,\moy{\mu}{\ABS{\GR g}^2}.
$$
Therefore, Poincar\'{e} inequality can be seen in a sense as the weakest
$\Phi$-Sobolev inequality since it is implied by any $\Phi$-Sobolev
inequality, at least on $\cC_b^\infty(\dR^d,\cI)$.  Moreover, let us assume
that $\dR_+\subset\cI$ and let us take now $f_\veps=a+\veps\,g^{p/2}$ where
$g$ is positive and $p\in(1,2]$, then we get by the same way:
$$
\moy{\mu}{g^p} - \moy{\mu}{g^{p/2}}^2%
\leq 2\,c\,p^2\,\moy{\mu}{\ABS{\GR g}^2\,g^{p-2}}.
$$
Now, since $p\in(1,2]$, Jensen inequality yields
$\moy{\mu}{g^{p/2}}^2\leq\moy{\mu}{g}^p$, and therefore:
$$
\moy{\mu}{g^p} - \moy{\mu}{g}^p%
\leq 2\,c\,\frac{p}{(p-1)}\,\moy{\mu}{p(p-1)\,g^{p-2}\,\ABS{\GR g}^2},
$$
which is a $\Phi$-Sobolev inequality with $\Phi(x)=x^p$ on $\cI=\dR_+$.
However, the constant does not give the sharp one obtained in Corollary
\ref{co:phi-sob-logconc-rd} when $\mu$ is $\log$-concave, and there is no way
to recover a logarithmic Sobolev inequality by letting $p$ tends to $1^+$ like
in \eqref{eq:xp-to-ent} since $2\,c\,p/(p-1)^2$ blows up near $p=1$. This fact
in not surprising since one can start from a Poincar\'{e} inequality in the
latter, which is known to be strictly weaker than logarithmic Sobolev
inequality (the simplest counter-example is given by the exponential
probability law on the real line, see also \cite[Chap. 6]{MR2002g:46132}).
Recall that the $\Phi$-Sobolev inequality obtained in Corollary
\ref{co:phi-sob-logconc-rd} for $\log$-concave probability measures include
the optimal Poincar\'{e} and logarithmic Sobolev inequalities by considering
the appropriate $\Phi$, namely:
\begin{equation}\label{eq:ils-logconc}
  \forall f \in \bW^{1,\infty}(\dR^d,\dR_+^*),\,%
  \ent{\mu}{f} %
  \leq\frac{1}{2\rho}\,\moy{\mu}{\frac{\ABS{\GR f}^2}{f}},
\end{equation}
and 
\begin{equation}\label{eq:ip-logconc}
  \forall f \in \bW^{1,\infty}(\dR^d,\dR),\,%
  \var{\mu}{f}
  \leq\frac{1}{\rho}\,\moy{\mu}{\ABS{\GR f}^2},
\end{equation}
and more generally, for any $p\in(1,2]$:
\begin{equation}\label{eq:ixp-logconc}
  \forall f \in \bW^{1,\infty}(\dR^d,\dR_+^*),\,
  \ENT{\mu}{x\mapsto x^p}{f}
  \leq\frac{1}{2\rho}\,\moy{\mu}{p(p-1)\,f^{p-2}\,\ABS{\GR f}^2}.
\end{equation}
Recall that $\entf{\mu}=\ENTF{\mu}{x\mapsto x\log x}$ and that
$\varf{\mu}:=\ENTF{\mu}{x\mapsto x^2}$.  The $\bW^{1,\infty}$-regularity is
not optimal, but we are not interested here in such aspects. Notice that
$f^{p-2}\,\ABS{\GR f}^2=4p^{-2}\,\ABS{\GR f^{p/2}}^2$. The particular constant
in $p$ appearing in \eqref{eq:ixp-logconc} allows to derive the logarithmic
Sobolev inequality by using the fact that for a positive $f$:
\begin{equation}\label{eq:xp-to-ent}
  \lim_{p\to1^+} \frac{\moy{\mu}{f^p}-\moy{\mu}{f}^p}{p-1} =\pd_{p=1}\ 
  \BRA{\moy{\mu}{f^p}-\moy{\mu}{f}^p} =\ent{\mu}{f}.
\end{equation}
Therefore, the logarithmic Sobolev inequality appears as a sort of limiting
case of a family of $\Phi$-Sobolev inequalities associated to $\Phi(x)=x^p$
with $p$ close to $1^+$, and the constant is sharp when $\mu$ is
$\log$-concave. In another direction, by an appropriate change of function,
one can reformulate the $p$-inequality \eqref{eq:ixp-logconc} as follows:
$$
\forall f \in \bW^{1,\infty}(\dR^d,\dR_+^*), 
\,\moy{\mu}{f^2}-\moy{\mu}{f^{2/p}}^p
\leq\frac{2(p-1)}{p\rho}\,\moy{\mu}{\ABS{\GR f}^2},
$$
which gives by denoting $q=2/p$, i.e. $q\in [1,2)$:
$$
\forall f \in \bW^{1,\infty}(\dR^d,\dR_+^*), 
\,\moy{\mu}{f^2}-\moy{\mu}{f^q}^{2/q}
\leq\frac{(2-q)}{\rho}\,\moy{\mu}{\ABS{\GR f}^2}.
$$
By removing the positivity condition on $f$ as explained in Remark
\ref{rm:non-neg-func}, we get finally that for any $q\in[1,2)$:
\begin{equation}\label{eq:ixq-latala}
  \forall f \in \bW^{1,\infty}(\dR^d,\dR),
  \, \moy{\mu}{f^2}-\moy{\mu}{\ABS{f}^q}^{2/q}
  \leq\frac{(2-q)}{\rho}\,\moy{\mu}{\ABS{\GR f}^2},
\end{equation}
which is exactly the inequality studied in \cite{MR89m:42027} for the Gaussian
measure and in \cite{MR2002b:60025}, \cite{wang-02}, \cite{MR1923685},
\cite{barthe-roberto} and \cite{barthe-cattiaux-roberto} for further
developments related to this particular case. Notice that in
\eqref{eq:ixq-latala}, the energy term $\moy{\mu}{\ABS{\GR f}^2}$ does not
depend on $q$, and one can adopt the following more convenient formulation:
\begin{equation}\label{eq:ixq-latala-sup}
  \forall\,f\in\bW^{1,\infty}(\dR^d,\dR),\,
  \sup_{q\in[1,2)}\,%
  \frac{\moy{\mu}{f^2}-\moy{\mu}{\ABS{f}^q}^{2/q}}{(2-q)}%
  \leq \frac{1}{\rho}\,\moy{\mu}{\ABS{\GR f}^2}.
\end{equation}
Generalisations of this type of statement where $(2-q)$ and $1/\rho$ are 
replaced by a more general functions are addressed in \cite{MR2002b:60025} and
\cite{wang-02}. Actually, inequality \eqref{eq:ixq-latala} appears as a sort
of infinite dimensional ``dual'' version of a Sobolev inequality.  Namely, let
$(\cM,g)$ be a smooth compact connected Riemannian manifold with dimension
$d\geq3$ and Ricci curvature $\rho>0$, and let $\mu$ be the normalised
Riemannian volume probability measure.  The Sobolev inequality states that for
any real valued smooth function $f$ on $M$:
\begin{equation}\label{eq:sob}
  \moy{\mu}{\ABS{f}^q}^{2/q} -\moy{\mu}{f^2} 
  \leq \frac{q+2}{2q}\,\frac{q-2}{\rho}\ \moy{\mu}{\ABS{\GR f}^2},
\end{equation}
where $q:=2d/(d-2)$ and $\ABS{\GR f}^2$ denotes the length of the gradient of
$f$ on $\cM$, see for example \cite{MR2002a:58045}.  Such an inequality is
stronger than logarithmic Sobolev inequality which is stronger than
Poincar\'{e} inequality. When $\cM$ is the standard $d$-dimensional sphere
$\dS^d(r)\subset\dR^{d+1}$ of radius $r=\sqrt{d}$, one has
$$
\rho=\frac{d-1}{r^2}=\frac{q+2}{2q}.
$$
It is then well known that one can recover the optimal logarithmic Sobolev
inequality for the standard Gaussian measure on $\dR^k$ when
$\mu=\cN(0,\rI_d)$ by taking the projection $\si_k:\dR^{d+1}\to\dR^k$ where
$k\leq d$ and by letting $d$ tends to $+\infty$ (i.e. $q$ tends to $2$).  This
fact is sometimes referred as the ``Poincar{\'e} observation'', see for
example \cite{MR2006070} for a generalisation to the case where $\mu$ is a
Boltzmann-Gibbs product measure.

\begin{rem}[From $\Phi$-Poincar\'{e} to $\Phi$-Sobolev]
  It is shown in \cite{wang-02} that for any $p\in[1,2]$ and any
  $f\in\rL^2(\Om,\mu,\dR)$ where $\mu$ is a probability measure on $\Om$ that:
  $$
  \moy{\mu}{f^2}-\moy{\mu}{\ABS{f}^p}^{2/p}%
  \leq\var{\mu}{f}+(1-p)\,\VAR{\mu}{x\mapsto{}x^p}{f}^{2/p}.
  $$
  It appears as an extended version of a Rothaus-Deuschel-Stroock-Bakry
  inequality to the case $p\in[1,2]$, cf. \cite[Lemmas 4.3.7 and
  4.3.8]{MR2002g:46132}. In presence of Poincar\'{e} inequality, it can be
  used perhaps to deduce a $\Phi$-Sobolev inequality for $\Phi(x)=x^p$ from an
  hypothetic $\Phi$-Poincar\'{e} inequality involving the $\Phi$-variance
  \eqref{eq:def-phi-var}.
\end{rem}

\begin{rem}[Sets convexity]
  For any $\al\in\{1,2,2'\}$, let $\dE(\al,\cI)$ be the set of smooth convex
  functions from $\cI$ to $\dR$ such that \HYPOT{$\al$} holds. It is easy to
  see that these three sets are convex vector cones, i.e. stable by linear
  combinations with non-negative scalar coefficients. For example, for any
  $\la:=(\la_1,\la_2,\la_3,\la_4)\in(\dR_+)^2\times\dR^2$ and any $p\in(1,2]$,
  the function $\Phi_{p,\,\la}:\dR_+\to\dR$ defined by
  $$
  \Phi_{p,\,\la}(x):=\la_1\,x^p+\la_2\,x\log x+\la_3\,x+\la_4
  $$
  is in $\dE(1,\dR_+)\,\cap\,\dE(2,\dR_+)\,\cap\,\dE(2',\dR_+)$. It is
  clear that functions $\Phi_{2,\,\la}$ are extremal points of the convex cone
  $\dE(1,\dR_+)$. One can observe that for any probability measure $\mu$ and
  any interval $\cI$ of $\dR$, the functional valued functional
  $\Phi\in\dE(\al,\cI)\mapsto\ENTF{\mu}{\Phi}$ is the restriction of a linear
  functional over the cone $\dE(\al,\cI)$, which is not sensitive to $\la_3$
  and $\la_4$. The same linearity property holds for the associated
  $\Phi$-Sobolev inequalities which are moreover invariant by any dilatation
  in $\la$. Thus, we understand that in view of $\Phi$-Sobolev inequalities,
  $\la_3$ and $\la_4$ plays no role and that $\Phi_{2,\,e_1}$ and
  $\Phi_{?,\,e_2}$ which corresponds respectively to classical variance and
  entropy are ``extremal'' cases under \HYPOT{1}. It could be nice to try to
  use a sort of Choquet integral representation by mean of extremal points,
  cf. \cite{MR1863703}.
\end{rem}

\section{Perturbation, tensorisation, convolution, concentration}
\label{se:pert-tenso-concent}

We explore here the stability of $\Phi$-Sobolev inequalities by tensorisation,
convolution, and perturbation.  We give also some concentration of measure
consequences.

\subsection{Tensorisation}
\label{ss:tensorisation}

A sub-additivity property for $\Phi$-entropies under \HYPOT{1} can be found in
\cite[Cor. 3]{MR2002b:60025}.  However, we provide in the sequel a short and
simple proof relying on the convexity of the $\Phi$-entropy functional
established in Proposition \ref{pr:convex-phi-ent}. Duality and sub-additivity
formulas are well known for variance and entropy, cf. \cite{MR97j:60005} and
\cite[Chap. 1]{MR2002g:46132}. The tensorisation property for $\Phi$-entropies
can be seen as an extension of the sub-additivity of Kullback-Leibler
relative entropy. Beware that it is different from the also well known
sub-additivity of Shannon like entropies, cf. Section \ref{se:shannon} page
\pageref{se:shannon} for an explanation.

Recall that for two probability measures $\nu$ and $\mu$ on a probability space
$(\Om,\cF)$ with $\nu\ll\mu$ , the Kullback-Leibler relative entropy 
$\ent{}{\nu\,\vert\,\mu}$ is defined by:
\begin{equation*}
  \ent{}{\nu\,\vert\,\mu}
  :=\int_\Om\frac{d\nu}{d\mu}\log\frac{d\nu}{d\mu}\,d\mu
  =\moy{\mu}{\frac{d\nu}{d\mu}\log\frac{d\nu}{d\mu}}
  =\moy{\nu}{\log\frac{d\nu}{d\mu}}.
\end{equation*}
Let $(\Om_1\times\cdots\times\Om_n, \cF_1\otimes\cdots\otimes\cF_n,
\mu_1\otimes\cdots\otimes\mu_n)$ be a product probability space. 
Then, for any probability measure $\nu \ll\mu_1\otimes\cdots\otimes\mu_n$:
\begin{equation}\label{eq:kullback-tenso}
  \ent{}{\nu\,\vert\,\mu_1\otimes\cdots\otimes\mu_n}
  \leq\moy{\mu}{\ent{}{\nu\,\vert\,\mu_1}}
  +\cdots+\moy{\mu}{\ent{}{\nu\,\vert\,\mu_n}},
\end{equation}
where 
$$
\ent{}{\nu\,\vert\,\mu_i} := \ent{\mu_i}{f}
= \int_{\Om_i}f\log f\,d\mu_i -\int_{\Om_i}\!f\,d\mu_i\log\int_{\Om_i}\!f\,d\mu_i,
$$
and $f:=d\nu/\mu_1\otimes\cdots\otimes\mu_n$. Moreover, equality holds if
and only if $\nu=\nu_1\otimes\cdots\otimes\nu_n$. In other words, for any real
valued integrable function $f$ on the product space:
\begin{equation}\label{eq:ent-tenso}
  \ent{\mu_1\otimes\cdots\otimes\mu_n}{f}
  \leq\moy{\mu}{\ent{\mu_1}{f}}+\cdots+\moy{\mu}{\ent{\mu_n}{f}}.
\end{equation}
This property of relative entropy remains true for $\Phi$-entropies when
$\Phi$ fulfils \HYPOT{1}:

\begin{prop}[Tensorisation of $\Phi$-entropies]
  \label{pr:phi-ent-tenso}
  Suppose that $\Phi$ satisfies \HYPOT{1}.  Let
  $(\Om_1\times\cdots\times\Om_n, \cF_1\otimes\cdots\otimes\cF_n,
  \mu_1\otimes\cdots\otimes\mu_n)$ be a product probability space.  Then, for
  any measurable function $f:(\Om,\cF)\to(\cI,\dB(\cI))$ one has:
  \begin{equation}\label{eq:phi-ent-tenso}
    \ENT{\mu_1\otimes\cdots\otimes\mu_n}{\Phi}{f} 
    \leq \moy{\mu}{\ENT{\mu_1}{\Phi}{f}} 
    +\cdots+\moy{\mu}{\ENT{\mu_n}{\Phi}{f}}.
  \end{equation}
\end{prop}
\begin{proof}
  As for variance and entropy, \eqref{eq:phi-ent-tenso} is a consequence of
  the convexity of the $\Phi$-entropy functional, and can be obtained by using
  the variational (duality) formula for $\Phi$-entropies
  \eqref{eq:phi-ent-varfor}.  Namely, for any $h$ and $i\in\{0,\ldots,n\}$,
  let $h_i:=\moy{\mu_1\otimes\cdots\otimes\mu_i}{h}$, with $h_0=h$. Thanks to
  the variational formula \eqref{eq:phi-ent-varfor}, one has:
  $$
  \ENT{\mu}{\Phi}{f}
  \leq
  \moy{\mu}{(\Phi'(h_0)-\Phi'(h_n))(f-h)}-\ENT{\mu}{\Phi}{h}.
  $$
  Since we can write
  $$
  \Phi'(h_0)-\Phi'(h_n) = \sum_{i=1}^n \PAR{\Phi'(h_{i-1})-\Phi'(h_i)},
  $$
  the desired result follows by the variational formula
  \eqref{eq:phi-ent-varfor} again and the fact that
  $$
  \ENT{\mu}{\Phi}{h} \geq \moy{\mu}{\ENT{\mu_i}{\Phi}{h_{i-1}}},
  $$
  which is due to Jensen inequality for the convex function $\Phi$.
\end{proof}

The tensorisation property for $\Phi$-entropies can be used to tensorise
$\Phi$-Sobolev inequalities as what is done for logarithmic Sobolev and
Poincar\'{e} inequalities, which can then be viewed as particular cases. In
essence, $\Phi$-Sobolev inequalities under \HYPOT{1} are then infinite
dimensional since they hold on the product space with the maximum of the one
dimensional constants.

%
\subsection{Convolution}
\label{ss:convolution}

For any $x\in\dR^d$, we denote by $\tau_x$ the translation of vector $x$ in
$\dR^d$ acting on any function $f:\dR^d\to\dR$ by
$(\tau_x\cdot{}f)(y):=f(x+y)$. One can state then the following immediate
Corollary of Proposition \ref{pr:phi-ent-tenso}.

\begin{cor}[Stability of $\Phi$-Sobolev inequalities by convolution]
  \label{co:phi-sob-convo}
  Suppose that $\Phi$ satisfies \HYPOT{1}.  For any $i\in\{1,\ldots,n\}$, let
  $\mu_i$ be a probability measure on $\dR^d$ satisfying to a $\Phi$-Sobolev
  inequality of the form
  $$
  \exists{}c_i>0,\quad\forall{}f:\dR^d\to\cI\text{ smooth},\quad%
  \ENT{\mu_i}{\Phi}{f} \leq c_i\,\moy{\mu_i}{\cE^\Phi(f)},
  $$
  where $f\mapsto\cE^\Phi(f)$ is a non-negative functional depending on
  $\Phi$ and its derivatives such that $\cE\circ\tau_x=\tau_x\circ\cE$ for any
  $x\in\dR^d$. Then, one has for any smooth function $f:\dR^d\to\cI$
  \begin{equation}
    \label{eq:phi-sob-convo}
    \ENT{\mu_1*\cdots*\mu_n}{\Phi}{f}%
    \leq (c_1+\cdots+c_n)\,\moy{\mu_1*\cdots*\mu_n}{\cE^\Phi(f)}.
  \end{equation}
\end{cor}
The commutation property $\cE\circ\tau_x=\tau_x\circ\cE$ is satisfied for
example when one has $\cE^\Phi(f)=\Phi''(f)\,\ABS{\GR{}f}^2$, and one can then
check the optimality of the obtained constant $c_1+\cdots+c_n$ when $\mu_i$
are Gaussian measures. The Poisson measures with suitable energy $\cE^\Phi$
give another example for which constant is optimal, cf. Section
\ref{se:phi-sob-levy}. We believe that Corollary \ref{co:phi-sob-convo}
remains essentially the same if one replaces the Euclidean space $\dR^d$ by an
infinite dimensional Banach space or by a topological Abelian group (and maybe
any Lie group). It can be useful for processes with independent increments, in
discrete or continuous time, and for infinitely divisible laws. However, we
will use different methods in the sequel, since i.i.d. increments and/or
infinite divisibility can be used by other ways to provide the same result.

\begin{rem}[Invariance by action of the translation group]
  \label{rm:inv-transl}
  Notice that hypothesis $\cE\circ\tau_x=\tau_x\circ\cE$ for any $x\in\dR^d$
  ensures that any associated $\Phi$-Sobolev inequality is ``invariant by
  translations in base space''. Namely, if one has
  $$
  \exists\,c>0,\quad\forall{}f\in\cA,\quad%
  \ENT{\mu}{\Phi}{f} \leq c\,\moy{\mu}{\cE^\Phi(f)},
  $$
  then for any $x\in\dR^d$, 
  $$
  \forall{}f\in\tau_{-x}\cdot\cA,\quad%
  \ENT{\tau_x\cdot\mu}{\Phi}{f} %
  \leq c\,\moy{\tau_x\cdot\mu}{\cE^{\Phi}(f)},
  $$
  where the action of $\tau_x$ on measure $\mu$ is given by
  $\tau_x\cdot\mu:=\mu*\de_x$.  Thus, if the class of functions $\cA$ is
  invariant by translations, i.e.  $\tau_x\cdot\cA=\cA$ for any $x\in\dR^d$,
  then any element of the orbit $\{\mu*\de_x,\,x\in\dR^d\}$ satisfies to the
  $\Phi$-Sobolev inequality satisfied by $\mu$, with same constant $c$ and
  class $\cA$.  Here again, we believe that things remain essentially the same
  if one replaces $\dR^d$ by an infinite dimensional Banach space or by a
  topological Abelian group.  Notice that when $\Phi(x)$ is of the form
  $x\log{}x$ or $x^p$ with $p\in(1,2]$, the associated $\Phi$-Sobolev
  inequality on $\dR^d$ with $\cE(f)=\Phi''(f)\,\ABS{\GR{}f}^2$ is
  homogeneous, and thus is additionally stable by dilatations.
\end{rem}

Actually, translations are particular examples of Lipschitz functions and 
one can expect a sort of stability by Lipschitz transforms when the energy
functional $\cE$ satisfies some stability. The following remark gives an
answer.

\begin{rem}[Invariance by action of Lipschitz functions]
  \label{rm:inv-lip}
  Let $\Te:\dR^d\to\dR^{d'}$ be a measurable map, acting on a function
  $f:\dR^{d'}\to\dR$ by $(\Te\cdot{}f)(y):=f(\Te(y))$ for any $y\in\dR^d$.
  Assume that the probability measure $\mu$ on $\dR^d$ satisfies to the
  following $\Phi$-Sobolev inequality:
  $$
  \exists\,c>0,\quad\forall\,g\in\cA,\quad%
  \ENT{\mu}{\Phi}{g} \leq c\,\moy{\mu}{\cE^\Phi(g)},
  $$
  where $\cE^\Phi$ is a non-negative functional. Assume that there exists a
  constant $\al>0$ such that for any $g\in\cA$,
  $(\cE^\Phi\circ\Te)(g)\leq{}\al\ (\Te\circ\cE^\Phi)(g)$. Let us denote by
  $\cA\cdot\Te$ the set of functions $f:\dR^{d'}\to\cI$ such that
  $\Te\cdot{}f\in\cA$.  Then, one has
  $$
  \forall\,f\in\cA\cdot\Te,\quad%
  \ENT{\Te\cdot\mu}{\Phi}{f} \leq c\,\al\,\moy{\Te\cdot\mu}{\cE^\Phi(f)},
  $$
  where $\Te\cdot\mu:=\mu\circ\Te^{-1}$ is the image measure of $\mu$ by
  $\Te$, i.e. the law of $\Te$ under $\mu$.  An important example is given by
  $\cE^\Phi(g)=\Phi''(g)\,\ABS{\GR{}g}^2$ and by any smooth map $\Te$ such
  that $\NRM{\mathrm{Jac(\Te)}}_2\leq{}\sqrt{\al}$. In particular, when
  $d=d'$, such $\Phi$-Sobolev inequalities are stable up to constants by the
  action of the non-Abelian group of diffeomorphisms with ``bounded
  Jacobian''.  Actually, in many cases including
  $\cE^\Phi(g)=\Phi''(g)\,\ABS{\GR{}g}^2$, only a weak smoothness of $\Te$ is
  needed, and for example one can assume only that $\Te$ is Lipschitz when
  $d'=1$.
\end{rem}

As an immediate consequence of Corollary \ref{co:phi-sob-convo} and Remarks
\ref{rm:inv-transl} and \ref{rm:inv-lip}, one can deduce the following result.

\begin{cor}\label{eq:phi-sob-stabil-gen}
  For any $i\in\{1,\ldots,n\}$, let $\mu_i$ be a probability measure on
  $\dR^{d_i}$ satisfying to a $\Phi$-Sobolev inequality of the form:
  $$
  \exists,c_i>0,\quad\forall{}f:\dR^{d_i}\to\cI\text{ smooth},\quad%
  \ENT{\mu_i}{\Phi}{f} \leq c_i\,\moy{\mu_i}{\Phi''(f)\,\ABS{\GR{}f}^2},
  $$
  where $\Phi$ satisfies \HYPOT{1}. Let $d$ be in $\dN^*$. For each
  $i=1,\ldots,n$, let $\Te_i:\dR^{d_i}\to\dR^d$ be a smooth (resp. Lipschitz
  when $d=1$) map with $\NRM{\mathrm{Jac}(\Te_i)}_2\leq\sqrt{\al_i}$ (resp.
  $\LIP{f}\leq\sqrt{\al_i}$ when $d'=1$). Let $\mu$ be the probability measure
  on $\dR^d$ defined by
  \begin{equation}\label{eq:kernel-mes}
    \mu:=(\Te_1\cdot\mu_1)*\cdots*(\Te_n\cdot\mu_n).
  \end{equation}
  Then, for any smooth function $f:\dR^{d}\to\cI$
  \begin{equation}\label{eq:phi-sob-kernel}
    \ENT{\mu}{\Phi}{f} \leq %
    (c_1\al_1+\cdots+c_n\al_n)\,\moy{\mu}{\Phi''(f)\,\ABS{\GR{}f}^2}.
  \end{equation}
\end{cor}

Probability measures like in \eqref{eq:kernel-mes} appear for example in
Statistics as fixed points of some kernel estimators. A very special Gaussian
case of \eqref{eq:kernel-mes} was studied in \cite{MR1914697} when
$\Phi(x)=x\log{}x$ and $(n,d_1,d)=(1,2,1)$. In view of some applications to
Statistics and Statistical Mechanics, the natural next step is to explore the
stability of $\Phi$-Sobolev inequalities by general mixtures of the form
$$
\moy{\mu}{f} := \int_{T}\,\moy{\mu_t}{f}\,d\nu(t);
$$
and conditioning of the form
$\mu=\cL((X_1,\ldots,X_n)\,\vert\,F(X_1,\ldots,X_n)=m)$ for example. But such
aspects are far more delicate and intricate and will be hopefully the aim of
forthcoming papers.

\subsection{Perturbation}
\label{ss:perturbation}

One can derive a perturbation property for the $\Phi$-Sobolev inequality via
the following straightforward variational formula for $\Phi$-entropies:
\begin{equation}\label{eq:phient-varfor-forperturb}
  \ENT{\mu}{\Phi}{f}=
  \inf_{a\in\inter{\cI}\subset\dR}\,
  \bE_\mu(\underbrace{\Phi(f)-\Phi(a)-\Phi'(a)\,(f-a)}_{\geq 0}).
\end{equation}
This formula is nothing else but the consequence of Taylor formula for the
convex function $\Phi$, and no more assumption on $\Phi$ are required here. It
can be found in \cite[Lem. 3.4.2]{MR2002g:46132}. One can easily recover the
well known formulas for variance and entropy by considering the appropriate
$\Phi$:
$$
\var{\mu}{f} = \inf_{a\in\dR}\moy{\mu}{(f-a)^2}
$$
and
$$
\ent{\mu}{f^2} = \inf_{a\in\dR_+^*}\moy{\mu}{f^2\log(f^2/a)-a+f^2}.
$$
Such a variational formula \eqref{eq:phient-varfor-forperturb} allows a
perturbation of $\Phi$-Sobolev inequalities via the method used by Holley and
Stroock for Poincar\'{e} and logarithmic Sobolev inequalities, cf.
\cite{MR89e:82013}. Namely, let $\mu$ be a probability measure on $(\Om,\cF)$
and $B:(\Om,\cF)\to(\dR,\dB(\dR))$ be a bounded measurable function. If we
define the probability measure $\nu_B$ on $(\Om,\cF)$ by:
\begin{equation}\label{eq:def-nub}
  d\nu_B := (Z_{\mu,B})^{-1}\,\exp(B)\,d\mu,
\end{equation}
where $Z_{\mu,B}:=\moy{\mu}{\exp(B)}$, then one can write:
\begin{align*}
\ENT{\nu_B}{\Phi}{f} 
&= \inf_{a\in\inter{\cI}} \,%
\int_\Om\!\!\PAR{\Phi(f)-\Phi(a)-\Phi'(a)\,(f-a)}\,d\nu_B\\
&\leq \exp\PAR{-\inf(B)+\sup(B)}\,\inf_{a\in\inter{\cI}}
\,\int_\Om\!\!\PAR{\Phi(f)-\Phi(a)-\Phi'(a)\,(f-a)}\,d\mu\\
&= \exp(\OSC{B})\,\ENT{\mu}{\Phi}{f}.
\end{align*}  
One can express the result as follows:
\begin{prop}[Perturbation]\label{th:phi-sob-pert}
  Let $(\Om,\cF,\mu)$ be a probability space such that
  $$
  \exists\,c\in\dR_+^*,
  \quad
  \forall f\in\cA,
  \quad
  \ENT{\mu}{\Phi}{f} \leq c\ \moy{\mu}{\cE(f)},
  $$
  where $\cA$ is a class of real valued measurable functions on $(\Om,\cF)$
  taking their values in $\cI$ and where $\cE:\cA\to\rL^1(\Om,\cF,\mu,\dR_+)$
  is a functional. Then, for any bounded measurable function
  $B:(\Om,\cF)\to\dR$, one has:
  $$
  \forall f\in\cA, \quad \ENT{\nu_B}{\Phi}{f} \leq c\ e^{2\,\OSC{B}}\ 
  \moy{\nu_B}{\cE(f)},
  $$
  where $\nu_B$ is like in \eqref{eq:def-nub}.
\end{prop}

\subsection{Concentration of measure}
\label{ss:concentration}

It is well known that a logarithmic Sobolev inequality for $\mu$ on $\dR^d$ of
the form
$$
\exists c>0,\quad\forall f:\dR^d\to\dR\text{ smooth},\quad%
\ENT{\mu}{x\mapsto x\log x}{f^2} \leq c\,\moy{\mu}{\ABS{\GR f}^2}
$$
gives, when applied to $f=\exp(\la F)$ where $F:\dR^d\to\dR$ is
$1$-Lipschitz, a Gaussian like exponential upper bound for the Laplace
transform of $\cL_\mu(F)$: $\moy{\mu}{\exp(\la F)} \leq \exp(c\,t^2)$.  Such a
bound can be then used via classical Cheby\-chev-Markov-Chernov-Cram\'{e}r
approach to give a concentration inequality for $F$ around its mean:
$$
\mu\PAR{\ABS{F-\moy{\mu}{F}}\geq t} \quad\leq\quad 2\,\exp(-t^2/c).
$$
The $t^2$ comes from the fact that the Young-Fenchel-Legendre convex
conjugate of $x\mapsto p^{-1}\,x^p$ is $x\mapsto q^{-1}\,x^q$ where
$q:=p/(p-1)$ is the H\"{o}lder conjugate of $p$, and thus $q=2$ when $p=2$.
This method is known as Herbst argument and gives precise and non-asymptotic
bounds which strengthen large deviations results. The concentration bound
obtained from logarithmic Sobolev in discrete settings are only Poissonian,
i.e. $\exp(-c\min(t^2,\,t\log t))$, due to the lack of chain rule. At the
opposite side of the set of possible $\Phi$ functions, the Poincar\'{e}
inequality gives by the same method an exponential like concentration, i.e.
$\exp(-c\,t)$.  See for example \cite{MR2002j:60002,MR1849347} and \cite[Chap.
7]{MR2002g:46132} for a general approach to concentration of measure via
functional inequalities.  It is tempting to study the concentration of measure
consequences of generic $\Phi$-Sobolev inequalities, and one can expect
intermediate exponential speeds between $t^2$ and $t$.  Let us recall what can
be found in \cite{MR2002b:60025} and \cite{MR1923685}. Let $\mu$ be a probability
measure on $\dR^d$, and $a\in [0,1]$ and $r=2/(2-a)$ (i.e. $r\in [1,2]$).
Assume that there exists a constant $C>0$ such that for any $q\in[1,2)$ and
any smooth function $f:\dR^d\to\dR$:
$$
\moy{\mu}{f^2}-\moy{\mu}{\ABS{f}^q}^{2/q}
\leq C\,(2-q)^a\,\moy{\mu}{\ABS{\GR f}^2}.
$$
We have seen already that such an inequality can be deduced by a change of
function from a $\Phi$-Sobolev inequalities with $\Phi(x)=x^{2/q}$. Then, the
following concentration of measure holds for any $t>0$ and $\mu$-integrable
$1$-Lipschitz function $F:\dR^d\to\dR$:
$$
\mu\PAR{F-\moy{\mu}{F} > \sqrt{C}\,t} \leq \exp\PAR{-K\,t^r},
$$
where $K>0$ is a universal constant. See
\cite{wang-02,MR1923685,MR2002b:60025,barthe-cattiaux-roberto} and references
therein for further developments. One can find a quite recent account in
\cite{lugosi-course}. Some aspects of concentration of measure consequences in
discrete space settings are addressed in \cite{MR2002f:60109} and
\cite{MR2002j:60002} and references therein.

Theorem \ref{th:phi-sob-diffusions} page \pageref{th:phi-sob-diffusions} gives
$\Phi$-Sobolev inequalities for hypercontractive diffusions and uniformly
strictly $\log$-concave probability measures, i.e. what is under Bakry-Emery
$\GD$ criterion \eqref{eq:curv-dim-crit}. It is quite natural to ask for a
general criterion to establish such inequalities beyond this scope. One can
find some answers in \cite{MR2002b:60025} and \cite{wang-02} for
$\Phi(x)=x^p$.  Namely, let $a\in [0,1]$ and $r=2/(2-a)$ (i.e. $r\in [1,2]$),
and consider the probability measure $\nu_r$ on $\dR^d$ defined by:
$$
d\nu_r(x):=Z_r^{-1}\,\exp\PAR{-\NRM{x}^r_r}\,dx.
$$
Then there exists an universal constant $C>0$ such that for any $q\in[1,2)$
and any smooth function $f:\dR^d\to\dR$:
$$
\moy{\nu_r}{f^2}-\moy{\nu_r}{\ABS{f}^q}^{2/q}
\leq C\,(2-q)^a\,\moy{\nu_r}{\ABS{\GR f}^2},
$$
which is exactly after the suitable change of functions the $\Phi$-Sobolev
inequalities for $\Phi(x)=x^{2/q}$.  One can try to investigate $\Phi$-Sobolev
inequalities on the real line via Hardy type inequalities, like what is known
for Poincar\'{e} and logarithmic Sobolev inequalities, cf. for example
\cite[Chap. 6]{MR2002g:46132}. As we will show in the sequel, $\Phi$-Sobolev
inequalities can be established on paths space of some diffusions (under
\HYPOT{1}) and some L\'{e}vy processes under \HYPOT{2}, extending by this way
what is already known in the literature for Poincar\'{e} and logarithmic
Sobolev inequalities.

\section{Phi-Sobolev for Brownian Motion and Wiener space}
\label{se:phi-sob-wiener}

Consider the standard Brownian Motion $(B_t)_{t\geq 0}$ on $\dR^d$ starting
from $B_0=0$. Since $\cL(B_t)=\cN(0,t\,\rI_d)$ is uniformly strictly
$\log$-concave with constant $1/t$, it satisfies a $\Phi$-Sobolev inequality
\eqref{eq:phi-sob-logconc} of constant $t/2$ under \HYPOT{1} for $\Phi$:
$$
\forall t>0,\ \forall f,\ 
\ENT{}{\Phi}{f(B_t)} \leq \frac{t}{2}\ \bE(\Phi''(f(B_t))\,\ABS{\GR f}^2(B_t)).
$$
Let $0 < t_1 < \cdots < t_n$ be $n$ successive times, and let
$F:\dR^n\to\dR$ be a smooth function. Then, one can write:
\begin{align*}
\ENT{}{\Phi}{F(B_{t_1},B_{t_2},\ldots,B_{t_n})} 
&= \ENT{}{\Phi}{F(Q_1,Q_1+Q_2,\ldots,Q_1+\cdots+Q_n)}\\
&= \ENT{}{\Phi}{G(Q_1,\ldots,Q_n)}\\
&= \ENT{\cL(Q_1,\ldots,Q_n)}{\Phi}{G},
\end{align*}
where $Q_i := B_{t_i}-B_{t_{i-1}}$ and $t_0:=0$ and
$$
G(x_1,\ldots,x_n):=F(x_1,x_1+x_2,\ldots,x_1+\cdots+x_n).
$$
But now, since Brownian Motion has i.i.d. increments, one gets
$$
\cL(Q_1,\ldots,Q_n)
=\cN(0,(t_1-t_0)\,\rI_d)\otimes\cdots\otimes\cN(0,(t_n-t_{n-1})\,\rI_d).
$$
Therefore, the tensorisation property \eqref{eq:phi-ent-tenso} yields the
following result.
\begin{thm}[Multi-times $\Phi$-Sobolev for Brownian Motion]
  Let $(B_t)_{t\geq 0}$ be a standard Brownian Motion on $\dR^d$.  Assume that
  $\Phi$ fulfils \HYPOT{1}. Then for any sequence of times
  $t_0:=0<t_1<\cdots<t_n$ and any smooth function $F:\dR^n\to\cI$:
 \begin{equation}\label{eq:phi-sob-multitime-bm}
   \ENT{}{\Phi}{F(B_{t_1},\ldots,B_{t_n})} %
   \leq \frac{1}{2}\,
   \moy{\cL(B_{t_1},\ldots,B_{t_n})}{\Phi''(F)\,\cD_{t_1,\ldots,t_n}^2 F},
 \end{equation}
 where
 \begin{equation}\label{eq:phi-sob-multitime-bm-diric}
 \cD_{t_1,\ldots,t_n}^2 F %
 := \sum_{i=1}^n (t_i-t_{i-1})\ \PAR{\sum_{j=i}^n \pd_j F}^2.
 \end{equation}
 Moreover, the inequality remains true when %
 $t_0:=0\leq t_1\leq \cdots\leq t_n$.
\end{thm}
One can ask if \eqref{eq:phi-sob-multitime-bm} is a consequence of
\eqref{eq:gaussian-brascamp-lieb}. If we define the $n\times n$ square
matrices $T$ and $Q$ by
$T:=\mathrm{diag}(\sqrt{t_1-t_0},\ldots,\sqrt{t_n-t_{n-1}})$ and
$Q_{i,j}:=\de_{j\geq{}i}$, then $\Si := (TQ)^\top TQ$ is the covariance matrix
of the centred Gaussian vector $(B_{t_1},\ldots,B_{t_n})$. In fact,
$\Si_{i,j}:=t_i\wedge t_j$, and the $\DP{\Si v;v}$ quadratic form in the right
hand side of \eqref{eq:gaussian-brascamp-lieb} reads
\begin{align}\label{eq:quad-form-bm-multi}
  \DP{\Si v;v} %
  &= (TQ v)^\top TQ v \nonumber\\
  &= \sum_{i=1}^n (t_{i}-t_{i-1})\,\PAR{\sum_{j=i}^n v_j}^2,
\end{align}
which is exactly the quadratic form appearing in the right hand side
\eqref{eq:phi-sob-multitime-bm-diric} of \eqref{eq:phi-sob-multitime-bm}.
Therefore, for Gaussian measures, the $\Phi$-Sobolev inequality
\eqref{eq:phi-sob-multitime-bm} obtained by applying the tensorisation formula
\eqref{eq:phi-ent-tenso} to the unidimensional case
\eqref{eq:phi-sob-logconc-Rn} (i.e. with $d=1$) is exactly the Brascamp-Lieb
type $\Phi$-Sobolev inequality \eqref{eq:gaussian-brascamp-lieb}. The direct
use of \eqref{eq:phi-sob-logconc-Rn} for the multivariate Gaussian measure
$\cL(B_{t_{1}},\ldots,B_{t_{n}})$ gives the quadratic form
$(\max\si(\Si))\,\sum_{i=1}^n v_i^2$ where $\si(\Si)$ is the spectrum of
$\Si=TQ$.  Thus, inequality \eqref{eq:phi-sob-logconc-Rn} for Gaussian
measures is a direct consequence of \eqref{eq:phi-sob-multitime-bm}. Notice
that $\det(\Si)=(t_1-t_0)\cdots(t_n-t_{n-1})$ and that
$\TR{\Si}=t_1+\cdots+t_n$. 

As we will see, inequality \eqref{eq:phi-sob-multitime-bm} appears as a
particular case of a more general one on paths space (i.e. for the Wiener
measure). The same procedure may be used for any random walk (resp. L\'{e}vy
process), provided that the law of the increment (resp. the infinitely
divisible law at time $t=1$) satisfies to a $\Phi$-Sobolev inequality under
\HYPOT{1} for $\Phi$ (we need a tensorisation property). Infinitely divisible
laws are particular cases of laws of sums of i.i.d. random variables. As we
will see, one can use for such laws the semi-group approach directly via the
associated L\'{e}vy process, to establish $\Phi$-Sobolev inequalities with
linear constant in time when $\Phi$ satisfies \HYPOT{2}, just like for
Brownian Motion.

One can derive the $\Phi$-Sobolev inequality for the Wiener measure on $\dR^d$
by a cylindrical method, starting from the multi-times $\Phi$-Sobolev
inequality \eqref{eq:phi-sob-multitime-bm} (established by tensorisation) by
letting the number $n$ of times considered tends to $+\infty$. However,
$\Phi$-Sobolev inequalities involving a Malliavin derivative can be easily
derived on paths space via a martingale representation approach and It\^{o}
formula as what was done for the $\log$-Sobolev inequality in
\cite{MR99b:60136}. Actually, one can state the following Theorem.
\begin{thm}[$\Phi$-Sobolev for the Wiener measure]%
  \label{th:phi-sob-wiener}
  Let $\bW_0(\dR^d)$ be the paths space of continuous
  functions from $[0,1]$ to $\dR^d$ starting from $0$ at $t=0$,
  measured by the standard Wiener measure.
  Assume that $\Phi$ fulfils \HYPOT{1}. Then for any random
  variable $F\in\rL^2(\bW_0(\dR^d),\,\cI)$:
  \begin{equation}\label{eq:phi-sob-wiener}
    \ENT{}{\Phi}{F} 
    \leq \frac{1}{2}\,\moy{}{\Phi''(F)\,\ABS{\cD F}_{\dH}^2},
  \end{equation}
  where $\cD:\rL^2(\mu)\to\rL^2(\mu,\dH)$ is the Malliavin gradient operator on
  $\bW_0(\dR^d)$ and $\dH$ is the Cameron-Martin Hilbert space.
\end{thm}
\begin{proof}
  For every functional $F : \bW_0(\dR^d)\to\cI$, consider the martingale
  $M_t:=\moy{}{F\,\vert\,\cF_t}$ where $0\leq t \leq 1$ and $\cF_t =
  \si(B_s,\,0\leq s\leq t)$ is the natural filtration. Then, one has
  $$
  d M_t = \DP{\moy{}{(\rD F)^\cdot_t\,\vert\,\cF_t},\,dB_t},
  $$
  where $(DF)_t^\cdot$ denotes the directional derivative of $F$. Now,
  It\^{o} formula yields:
  $$
  \moy{}{\Phi(M_t)}-\moy{}{\Phi(M_0)} 
  =\frac{1}{2}\, \moy{}{\int_0^1\!\! \Phi''(M_t)
    \,\ABS{\moy{}{(\rD F)^\cdot_t\,\vert\,\cF_t}}^2\,dt},
  $$
  which is nothing else but
  $$
  \ENT{\cL\PAR{(B_t)_{0\leq t\leq 1}}}{\Phi}{F} = \frac{1}{2}\,
  \moy{}{\int_0^1\!\!\Phi''\PAR{\moy{}{F\,\vert\,\cF_t}} 
    \,\ABS{\moy{}{(\rD F)^\cdot_t\,\vert\,\cF_t}}^2\,dt}.
  $$
  Now, by Jensen inequality for the bivariate convex function
  $(u,v)\mapsto\Phi(u)''\,v^2$, we get as in the semi-group proof of Theorem
  \ref{th:phi-sob-diffusions}:
  $$
  \Phi''\PAR{\moy{}{F\,\vert\,\cF_t}} 
  \,\ABS{\moy{}{(\rD F)^\cdot_t\,\vert\,\cF_t}}^2
  \leq
  \moy{}{\Phi''(F)\,\ABS{(\rD F)^\cdot_t}^2\,\vert\,\cF_t},
  $$
  which gives the desired result. Here again, one can use alternatively the
  concavity of $1/\Phi''$ and Cauchy-Schwarz inequality.
\end{proof}
Notice that the proof is a replica at paths space level of the semi-group
proof for diffusions in Theorem \ref{th:phi-sob-diffusions}, where $\SGf{t}$
is ``replaced'' by $\moy{}{F\,\vert\,\cF_t}$ and the diffusion property by
It\^{o} formula. This analogy is not formal indeed, since the diffusion
semi-group gives the solution of Stroock-Varadhan martingale problem
associated to the related elliptic diffusion operator, see for example
\cite[Chap. 5, Sect. 4]{MR92h:60127}. Such a proof can be extended to paths
spaces on manifolds as what was already done for the logarithmic Sobolev
inequality in \cite{MR99b:60136,MR98j:60110} (see also \cite[Sect.
8.3]{MR1882015} and references therein), as stated in the following Theorem.

\begin{thm}[$\Phi$-Sobolev for Brownian Motion on a Manifold]
  \label{th:phi-sob-bm-manif}
  Let $(\cM,g)$ be a smooth complete and connected Riemannian manifold
  equipped with the Levi-Civita connection. Let $x\in\cM$ and $W_{x}(\cM)$ be
  the space of continuous paths $\ga:[0,1]\to\cM$ with $\ga(0)=x$. Assume that
  the Ricci curvature of $\cM$ is uniformly bounded by the real number $K$.
  Assume that $\Phi$ fulfils \HYPOT{1}, then for any smooth non-negative
  random variable $F$ on $W_{x}(\cM)$:
  \begin{equation}\label{eq:phi-sob-bm-manif}
    \ENT{}{\Phi}{F} \leq%
    \frac{1}{2}\,e^K\,\moy{}{\Phi''(F)\,\ABS{\cD F}_{\dH}^2}.
  \end{equation}
\end{thm}
When $\cM$ is the standard $d$-dimensional Euclidean space $\dR^d$, one has
$K=0$ and we recover Theorem \ref{th:phi-sob-wiener}. We believe that
$\Phi$-Sobolev inequalities under \HYPOT{1} still hold on the paths space of
diffusions on manifolds with Driver total antisymmetry condition, as what
was done in \cite{MR99b:60136} for the logarithmic Sobolev inequality. We have
to mention that this beautiful martingale method for Wiener measure over
Riemannian manifold appeared for the first time in \cite{MR94m:58238} for the
Poincar\'{e} inequality. We will use roughly the same method under \HYPOT{2}
on Poisson space in Section \ref{ss:phi-sob-ps} page \pageref{ss:phi-sob-ps}.

\section{Phi-Sobolev for pure jumps L\'{e}vy processes and Poisson space}
\label{se:phi-sob-levy}

It is tempting to try to establish $\Phi$-Sobolev inequalities on Poisson
space. In particular, L\'{e}vy processes have i.i.d. increments like Brownian
Motion. Such processes are not diffusions, and the lack of the chain rule
\eqref{eq:chain-rule} forbids the direct use of the Bakry-Emery semi-group
proof. More precisely, for all Markov processes, $\GA \geq 0$, and the
$\mathrm{CD}(\rho,\infty)$ criterion is equivalent to $\GD \geq \rho\,\GA$.
But the equivalent form in term of strong commutation
\eqref{eq:strong-commut-sg} between the semi-group and $\GA$ is available only
for diffusions. Nevertheless, for L\'{e}vy processes, one has $\GD\geq 0$ and
indeed $\GD\geq\GA\sqrt{\GA}$, which gives Brownian Motion like commutation
formulas and then permits to derive respectively Poincar\'{e} and
$\log$-Sobolev inequalities, see \cite{MR2001i:60094} for the simple Poisson
point process and \cite{chafai-malrieu-levy} for more general L\'{e}vy
processes. More general approaches are presented in \cite{MR2002f:60109}
(generic Poisson space) and \cite{MR2002h:60087} (normal martingales), but
here again only for Poincar\'{e} and $\log$-Sobolev inequalities.

Let us start with the simplest pure jumps L\'{e}vy process $(X_t)_{t\geq 0}$
which is the simple Poisson point process with intensity $\la>0$, for which
\begin{equation}\label{eq:def-gi-levy}
(\GI f)(x):=\la\,\PAR{f(x+1)-f(x)}=:\la\,(\rD_1 f)(x),
\end{equation}
where $(\rD_y f)(x):=f(x+y)-f(x)$. If $g:=\SG{t-s}{f}$, one has:
$$
\pd_s\,\SG{s}{\Phi(g)} := \SG{s}{\GI \Phi(g)-\Phi'(g)\,\GI g}
=\la\,\SG{s}{\Psi(g,\,\rD_1 g)},
$$
where $\Psi$ is defined by \eqref{eq:def-psi}. But now, under \HYPOT{2},
$\Psi$ is bivariate convex on $\DBLI$. It was already observed for the simple
cases $\Phi(x)=x^2$ and $\Phi(x)=x\log x$ in the pretty paper
\cite{MR2002f:60109}. Hence, by the bivariate Jensen inequality and the
commutativity property between $\rD_1$ and $\SGf{t-s}$ (due to the i.i.d.
nature of the increments), one gets the following result.
\begin{thm}[Local $\Phi$-Sobolev for the simple Poisson point process]
  \label{th:phi-sob-loc-sppp}
  Let $(X_t)_{t\geq 0}$ be the simple Poisson point process on $\dR^d$ with
  intensity $\la>0$.  Assume that $\Phi$ fulfils \HYPOT{2}, then for any $t>0$
  and any smooth $f:\dR^d\to\cI$:
  \begin{equation}\label{eq:phi-sob-loc-sppp}
    \ENT{\SGf{t}}{\Phi}{f} \leq \la\,t\,\SG{t}{\Psi(f,\,\rD_1 f)},
  \end{equation}
  where $\SG{t}{f}(x):=\bE(f(X_t)\,\vert\,X_0=x)$ is the associated semi-group
  as in \eqref{eq:def-sg}.
\end{thm}
In particular, by taking $t=1$, we get:
\begin{cor}[$\Phi$-Sobolev for the simple Poisson measure]\label{co:phi-sob-mes-pois}
  Let $\cP_\la$ be the Poisson measure $\cP_\la$ of mean $\la>0$, then, under
  \HYPOT{2}:
  \begin{equation}\label{eq:phi-sob-mes-pois}
    \ENT{\cP_\la}{\Phi}{f}
    \leq \la\,\moy{\cP_\la}{\Psi(f,\,\rD_1 f)}.
  \end{equation}
\end{cor}
Such a result remains true for pure jumps L\'{e}vy processes (i.e.  without
Brownian part) by replacing $\rD_1$ by the appropriate jump integral. The
Brownian part may be added at final stage via tensorisation when possible.
Namely, we can state the following.
\begin{thm}[Local $\Phi$-Sobolev for a pure jumps L\'{e}vy process]
  \label{th:phi-sob-loc-levy}
  Let $(X_t)_{t\geq 0}$ be a pure jump L\'{e}vy process with infinitesimal
  generator of the form:
  \begin{equation}\label{eq:gi-loc-levy}
    (\GI f)(x) := \la\,\int_{\dR^d}\!\!
    \SBRA{\rD_y{f}(x)-\te\,\frac{y}{1+\ABS{y}^2}\cdot\GR f}\,d\nu(y), 
  \end{equation}
  where $(\la,\te)\in\dR_+^*\times\dR^d$ and where $\nu$ is a L\'{e}vy measure
  on $\dR^d$.  Assume that $\Phi$ fulfils \HYPOT{2}. Then for any $t>0$ and
  any smooth $f:\dR^d\to\cI$:
  \begin{equation}\label{eq:phi-sob-loc-levy}
    \ENT{\SGf{t}}{\Phi}{f}%
    \leq \la\,t\, \SG{t}{\int_{\dR^d}\!\!\Psi(f,\,\rD_y{}f)\,d\nu(y)},
  \end{equation}
  where $\SG{t}{f}(x):=\bE(f(X_t)\,\vert\,X_0=x)$ is the associated
  semi-group as in \eqref{eq:def-sg}.
\end{thm}
Taking $t=1$ in \eqref{eq:phi-sob-loc-levy} gives the same functional
inequality for the infinitely divisible law $\cL(X_1)$, exactly like in
Corollary \ref{co:phi-sob-mes-pois} for the simple Poisson measure. Theorem
\ref{th:phi-sob-loc-sppp} for the simple Poisson point process can be
recovered by taking $(\te,\nu)=(0,\de_1)$ in Theorem
\ref{th:phi-sob-loc-levy}. 

One can state the following multi-times version of Theorem
\ref{th:phi-sob-loc-levy} via the so called Lu-Yau-Stroock-Zegarli{\'n}ski
Markov decomposition method.
\begin{thm}[Multi-times $\Phi$-Sobolev for pure jumps L\'{e}vy proc.]%
  \label{th:phi-sob-multitime-levy}
  Let $(X_t)_{t\geq 0}$ be a pure jumps L\'{e}vy process on $\dR^d$ as in
  Theorem \ref{th:phi-sob-loc-levy}. Assume that $\Phi$ fulfils \HYPOT{2}.
  Then for any increasing sequence of times $t_0:=0<t_1<\cdots<t_n$ and any
  smooth $F:\dR^n\to\cI$:
  \begin{equation}\label{eq:phi-sob-multitime-levy}
    \ENT{}{\Phi}{F(X_{t_1},\ldots,X_{t_n})} \leq \la\, %
    \moy{}{\cD^\Phi\PAR{F(X_{t_1},\ldots,X_{t_n})}},
  \end{equation}
  where
  $$
  \cD^\Phi(F) %
  := \int_{\dR^d}\sum_{i=1}^n(t_i-t_{i-1})\,%
  \Psi(F,\rD_y^{i,\ldots,n} F)\,d\nu(y),
  $$
  where for any $x\in\dR^n$
  $$
  \rD_y^{i,\ldots,n} F(x):= F\circ\tau_i(y)(x)-F(x),
  $$
  and where for any $i\in\{1,\ldots,n\}$ and $x,y\in\dR^d$
  $$
  F\circ\tau_i(y)(x):=F(x_1,\ldots,x_{i-1},x_i+y,\ldots,x_n+y).
  $$
\end{thm}
\begin{proof}
  By induction on $n$, we can restrict the problem to the case $n=2$. Let $s$
  and $t$ be two distinct times with $s<t$. We assume that for any $u>0$, a
  $\Phi$-Sobolev inequality holds for $\cL(X_u)$ (i.e. for $\SGf{u}$) with
  constant $\la u$.  We would like to obtain a similar inequality for
  $\cL(X_s,X_t)$.  Let $F:\dR^2\to\cI$ be a smooth bivariate functional.  We
  start with the following conditional decomposition of the $\Phi$-entropy
  (which replaces in some ways the tensorisation property):
  \begin{align*}
    \ENT{}{\Phi}{F(X_s,X_t)} 
    =&\bE\SBRA{\bE(\Phi(F)\,\vert\,X_s)
      -\Phi(\bE(F\,\vert\,X_s))} \\
    &+\bE(\Phi(\bE(F\,\vert\,X_s)))
    -\Phi(\bE(\bE(F\,\vert\,X_s))),
  \end{align*}
  where we abridged $F(X_s,X_t)$ in $F$ in the right hand side. In other
  words,
  \begin{equation}\label{eq:decomp-phi-ent}
    \ENT{}{\Phi}{X}%
    =\bE(\ENT{}{\Phi}{X\,\vert\,Y)}+\ENT{}{\Phi}{\bE(X\,\vert\,Y)},
  \end{equation}
  where $\ENT{}{\Phi}{X\,\vert\,Y}$ is the conditional $\Phi$-entropy defined
  by
  $$
  \ENT{}{\Phi}{X\,\vert\,Y}:=\bE(\Phi(X)\,\vert\,Y)-\Phi(\bE(X\,\vert\,Y)).
  $$
  Now, if $G_x(y):=F(x,x+y)$, we have for any $x$:
  $$
  \ENT{}{\Phi}{G_x(X_{t-s})} \leq \la\,(t-s)\,
  \bE\PAR{\int_{\dR^d}\!\!\Psi(G_x(X_{t-s}),\rD_z G_x(X_{t-s}))\,d\nu(z)}.
  $$
  Therefore, since $X_t=X_s+X_t-X_s$ and
  $\cL(X_s,X_t-X_s)=\cL(X_s)\otimes\cL(X_{t-s})$, the first term of the
  $\ENTF{}{\Phi}$ decomposition can be bounded above as
  $$
  \bE\SBRA{\bE(\Phi(F)\,\vert\,X_s)-\Phi(\bE(F\,\vert\,X_s))} \leq
  \la\,(t-s)\,\bE\PAR{\int_{\dR^d}\!\!\Psi(F,\rD_z^2 F)\,d\nu(z)},
  $$
  where the $2$ exponent in $\rD_z^2$ means that the $z$ translation is
  done on the second variable of $F$ in the definition of $\rD_z$. On the
  other hand, if we define $H_{t-s}$ by
  $$
  H_{t-s}(x):=\bE(F(x,x+X_{t-s})),
  $$
  we have:
  $$
  \ENT{}{\Phi}{H_{t-s}(X_s)}
  \leq \la\,s
  \,\bE\PAR{\int_{\dR^d}\!\!\Psi(H_{t-s}(X_s),\rD_z H_{t-s}(X_s))\,d\nu(z)}.
  $$
  Now, by the commutation formula $\rD_z H_{t-s}(x):=\bE(\rD_z (x\mapsto
  F(x,x+X_{t-s})))$ and the bivariate Jensen inequality, the last term of the
  $\ENTF{}{\Phi}$ decomposition can be bounded above as follows:
  $$ 
  \bE(\Phi(\bE(F\,\vert\,X_s)))-\Phi(\bE(\bE(F\,\vert\,X_s)))
  \leq \la\,s\,\bE\PAR{\int_{\dR^d}\!\!\Psi(F,\rD_z^{1,2} F)\,d\nu(z)},
  $$
  where this time the translation is done in both variables.
\end{proof}
In some sense, the Markovian decomposition used in the proof of Theorem
\ref{th:phi-sob-multitime-levy} replaces the tensorisation property
\eqref{eq:phi-ent-tenso} when \HYPOT{2} holds instead of \HYPOT{1}. It was
introduced in Statistical Mechanics for Poincar\'e and logarithmic Sobolev
inequalities for finite volume Boltzmann-Gibbs measures related to spins
systems, cf. for example \cite{guionnet-zegarlinski-ihp-2002} and
\cite{MR2000k:60002}.

As for Brownian Motion, one can use a cylindrical method letting $n$ tends to
$+\infty$ in \eqref{eq:phi-sob-multitime-levy} in order to obtain a
$\Phi$-Sobolev inequality on paths space, as expressed in the following
Theorem. 
\begin{thm}[$\Phi$-Sobolev on paths space of pure jumps L\'{e}vy proc.]%
  \label{th:phi-sob-paths-levy}
  Let $(X_t)_{t\geq 0}$ be a pure jump L\'{e}vy process as in Theorem
  \ref{th:phi-sob-loc-levy}.  Assume that $\Phi$ fulfils \HYPOT{2}. Then, for
  any suitable function $F$ of $(X_t)_{t\geq 0}$ taking its values in $\cI$
  and any $T>0$:
  \begin{equation}\label{eq:phi-sob-paths-levy}
    \ENT{\cL\PAR{\PAR{X_t}_{0\leq t \leq T}}}{\Phi}{F}
    \leq \la\,\moy{}{\int_0^T\!\!\!\int_{\dR^d}\Psi(F,\rD_y^t F)\,d\nu(y)\,dt},
  \end{equation}
  where
  $$
  \rD_y^t F ((x_s)_{0\leq s\leq T})
  := F((x_s+y\,\rI_{[t,T]}(s))_{0\leq s\leq T})-F((x_s)_{0\leq s\leq T}).
  $$
\end{thm}
Obviously, Theorem \ref{th:phi-sob-paths-levy} is stronger than Theorem
\ref{th:phi-sob-loc-levy} since the latter can be recovered by simply
considering in \eqref{eq:phi-sob-paths-levy} functions $F$ of the form
$$
F((x)_{0\leq{}s\leq{}T})=G(x_{t_1},\ldots,x_{t_n})
$$
where $T$ is chosen bigger than $t_n$.  Inequality
\eqref{eq:phi-sob-multitime-levy} is similar to the one obtained for Brownian
Motion \eqref{eq:phi-sob-multitime-bm} via tensorisation. Actually, one can
use Lu-Yau-Stroock-Zegarli{\'n}ski Markovian method for Brownian Motion.
Generic L\'{e}vy processes -- and thus Brownian Motion -- are particular
examples of normal martingales and we believe that the approach used in
\cite{MR2002h:60087} for logarithmic Sobolev inequalities remains valid for
$\Phi$-Sobolev inequalities, but one has to precise the condition on $\Phi$.

As explained below, $\Phi$-Sobolev inequalities under \HYPOT{2} can be
established on paths space for generic Poisson space via a
Clark-Oc\^{o}ne-Haussmann formula.

\subsection{Phi-Sobolev inequalities on Poisson space}\label{ss:phi-sob-ps}
%


Let us explain finally how one can recover the $\Phi$-Sobolev inequality on
Poisson space via a Clark-Oc\^{o}ne-Haussmann formula like what is done in
\cite{MR2002f:60109} for the logarithmic Sobolev and Poincar\'{e} inequalities
(see also \cite{MR2002h:60087}). Following closely \cite{MR2002f:60109}, let
$(E,\cB,\nu)$ be a measured space where $\nu$ is a $\si$-finite measure and
$(W,\cF,\cP)$ the associated Poisson space with compensation measure $\nu$.
Let also $(\Om,\cC, \bP)$ be the Poisson space associated to $([0,1]\times
E,\cB([0,1])\otimes \cB, \mu(dt,dz):=dt\otimes d\nu(z))$, with compensation
measure $\mu$. One can then define the difference operators
$\rD_z:\rL^0(W,\cP)\to\rL(E\times W,\nu\otimes\cP)$ and
$\rD_{t,z}:\rL^0([0,1]\times E,\mu)\to\rL^0([0,1]\times E \times \Om,
\mu\otimes \bP)$ by:
$$
\rD_z F(\om) := F(\om+\de_z)-F(\om)
\quad\text{and}\quad
\rD_{t,z} \widehat{F}(\om) := \widehat{\rD_z F}(\om),
$$
where $\widehat{F}(\om):=F(\om([0,1],dz))$.  For any $t\in[0,1]$, we define
$$
\cC_t:=\si\PAR{\om(A); A\in\cB([0,t])\otimes \cB}.
$$
Let now $G\in\rL^2(\Om,\bP)$ and let $g(t,z,\om)$ be a predictable
$dt\otimes d\nu(z)\otimes\bP(d\om)$ version of
$\moy{\bP}{\rD_{t,z}\,G\,\vert\,\cC_t}$. One has then the following
Clark-Oc\^{o}ne type predictable representation formula :
\begin{equation}\label{eq:clark-ocone}
  G=\moy{\bP}{G}+\int_0^1\!\!\!\int_E\!g(t,z,\cdot)\,d\tilde{\om}(t,z),
\end{equation}
where $\tilde{\om}=\om-\mu$ is the compensated Poisson point process and
where the integral is taken in the sense of It\^{o}. Let us assume for
convenience that $G>\inf(\cI)$ and let us define the right continuous
martingale $(M_t)_{t\in[0,1]}$ by $M_t:=\moy{\bP}{G\,\vert\,\cC_t}$. Now,
since $M_{t^-}=M_t$ $dt\otimes\bP$-a.s., It\^{o} formula for jumps processes
(cf. \cite{MR90m:60069}) gives:
\begin{align*}
\ENT{\bP}{\Phi}{G}%
&=\moy{\bP}{\Phi(M_1)-\Phi(M_0)}\\%
&=\moy{\bP}{\int_0^1\!\!\int_E\!\Psi(M_t,g(t,z))\,dt\,d\nu(z)}.
\end{align*}
Now, it remains to use the bivariate convexity of $\Psi$ which comes from
\HYPOT{2} in order to get via Jensen inequality and \eqref{eq:clark-ocone}
that:
$$
\Psi(M_t,g(t,z))%
= \Psi\PAR{\moy{\bP}{G\,\vert\,\cC_t},\moy{\bP}{\rD_{t,z}\,G\,\vert\,\cC_t}}%
\leq \moy{\bP}{\Psi(G,\rD_{t,z}\,G)\,\vert\,\cC_t},
$$
which gives finally: 
\begin{equation}\label{eq:phi-sob-gps}
\ENT{\bP}{\Phi}{G}%
\leq \moy{\bP}{\int_0^1\!\!\!\int_E\!\Psi(G,\rD_{t,z}\,G)\,dt\,d\nu(z)}.
\end{equation}
The time interval $[0,1]$ can be easily replaced by $[0,T]$. One can recover
Theorems \ref{th:phi-sob-loc-sppp}, \ref{th:phi-sob-multitime-levy},
\ref{th:phi-sob-paths-levy} and \ref{th:phi-sob-multitime-levy} as
Corollaries, in exactly the same way used in \cite{MR2002f:60109} for
Poincar\'{e} and logarithmic Sobolev inequalities. Additionally, one can
derive by the same method F.K.G. inequalities, as what is done in
\cite{MR2002f:60109} for Poisson space and in \cite{MR2001c:60124} for
diffusions. As we have seen, the method used here for Poisson space is a
replica of the method used to establish Theorem \ref{th:phi-sob-wiener}
concerning Wiener measure. The major difference is the lack of chain rule in
discrete space settings which leads to replace \HYPOT{1} by \HYPOT{2}.

\subsection{Some remarks about discrete Dirichlet forms}\label{ss:rem-dirfors-discr}

We collect here few remarks about \HYPOT{2}, \HYPOT{2'}, and comparisons with
standard Dirichlet forms in discrete space settings.

\begin{rem}[Comparison with classical Dirichlet forms]\label{rm:dircomp}
  Assume that \HYPOT{2'} holds.  Since $\Phi'$ is concave, $\Phi''$ is
  non-increasing and therefore, when $v\geq 0$
  $$
  \Psi(u,v) 
  = \frac{1}{2}\,\int_u^{u+v}\!\!(u+v-w)\,\Phi''(w)\,dw \leq \Phi''(u)\,v^2.
  $$
  The case $v\leq 0$ is very similar, and we get finally that under \HYPOT{2'}:
  \begin{equation}\label{eq:dircomp-df-df}
    \Psi(u,v) \leq \Phi''(u)\,v^2 \text{\quad on \quad} \DBLI.
  \end{equation}
  On the other hand, if we assume only \HYPOT{2}, one can write:
  $$
  \Psi(u,v) %
  = \frac{1}{2}\,\int_u^{u+v}\!\!(u+v-w)\,\Phi''(w)\,dw %
  \leq v\,\PAR{\Phi'(u+v)-\Phi'(u)},
  $$
  which gives that under \HYPOT{2}:
  \begin{equation}\label{eq:dircomp-df-dlogf}
    \Psi(u,v) \leq v\,\PAR{\Phi'(u+v)-\Phi'(u)} \text{\quad on \quad} \DBLI.
  \end{equation}
  Thus, under \HYPOT{2} (resp. \HYPOT{2'}), we recover for example the well
  known $(\rD f)^2/f$ (resp. $\rD f\,\rD\log f$) Dirichlet forms when
  $\Phi(x)=x\log x$, as in \cite{MR2002f:60109}. When $\Phi(x)=x^2$, both of
  them give $2\,(\rD f)^2$. 
\end{rem}

\begin{rem}[``L-1'' and ``L-2'' forms of logarithmic Sobolev inequalities]
  In continuous settings, the equivalence between $\rL^1$ and $\rL^2$ forms of
  the logarithmic Sobolev inequality, given respectively by
  \begin{equation}\label{eq:lsi-l2}
  \forall f,\ \ent{\mu}{f^2} \leq c\,\moy{\mu}{\GA f}
  \end{equation}
  and
  \begin{equation}\label{eq:lsi-l1}
  \forall f>0,\ \ent{\mu}{f} \leq \frac{1}{4}\,c\,\moy{\mu}{\frac{\GA f}{f}}  
  \end{equation}
  is a consequence of the chain rule (\ref{eq:chain-rule}) for $\GA$, which is
  itself a consequence of the diffusion property \eqref{eq:def-diffusion} for
  the associated infinitesimal generator. In discrete settings, the lack of
  chain rule destroys this equivalence, and actually, the simple Poisson
  measure satisfies the $\rL^1$ form but not the $\rL^2$ form, cf. for example
  \cite[Chap.  1]{MR2002g:46132}. The fact that the semi-group of the simple
  Poisson point process is not hypercontractive was noticed by Surgailis
  eighteen years ago in \cite{MR86k:60098}.  The concentration of measure
  consequences of the $L^1$ form are only Poisson-like in discrete settings,
  which is completely logical, cf. \cite{MR2002j:60002}.
\end{rem}

\begin{rem}[Poissonian $\rL^1$ logarithmic Sobolev inequality]
  \label{rm:dircompls}  
  Notice that the $\GA$ operator (cf.  \eqref{eq:def-gamma}) associated to the
  pure jumps L\'{e}vy process with generator $\GI$ given by
  \eqref{eq:def-gi-levy} is:
  \begin{equation}\label{eq:gamma-discr}
  (\GA f)(x):=\frac{\la}{2}\,\int_{\dR^d}\!\! \ABS{\rD_y f(x)}^2\,d\nu(y).
  \end{equation}
  When $\Phi(x)=x\log x$, the Dirichlet forms comparison
  \eqref{eq:dircomp-df-dlogf} yields
  \begin{align*}
    \ent{\SGf{t}}{f}%
    &\leq\la t\,\SG{t}{\int_{\dR^d}\!\!f^{-1}(\rD_y f)^2\,d\nu(y)} \\
    &=\la t\,\SG{t}{f^{-1}\int_{\dR^d}\!\!(\rD_y f)^2\,d\nu(y)} \\
    &=2t\,\SG{t}{\frac{\GA f}{f}},
  \end{align*}
  which is the $\rL^1$ form \eqref{eq:lsi-l2} for measure $\SGf{t}$ and with
  constant $2t$.
\end{rem}

\begin{rem}[Discussion on \HYPOT{2} and \HYPOT{2'}]\label{rm:h1-h2}
  Let $\zeta_1(u,v)=\Phi''(u)\,v^2$ and
  $\zeta_2(u,v)=(\Phi'(u+v)-\Phi'(u))\,v$, for $(u,v)\in\DBLI$.  Such
  functions are non negative since $\Phi$ is convex.  One has $\Psi \leq
  \zeta_1$ when $\Phi'$ is concave, and since
  $$
  \GR^2\!\zeta_1(u,v) = 
  \PAR{
    \begin{array}{rcl}
      \Phi''''(u)\,v^2   & \quad & 2\,\Phi'''(u)\,v \\
      2\,\Phi'''(u)\,v   & \quad & 2\,\Phi''(u)  
    \end{array}},
  $$
  \HYPOT{1} is equivalent to $\Phi''''\,\Phi'' \geq 2\,\Phi'''^2$ which is
  nothing else than the concavity of $1/\Phi''$ when $\Phi''$ is convex. One
  can check that $\Phi(x)=x\log x$ and $\Phi(x)=x^2$ are in some sense
  extremal solutions of this O.D.I., up to affine additions. Similarly, one
  has $\Psi \leq \zeta_2$ when $\Phi$ is convex, and since
  $\GR^2\!\zeta_2(u,v)$ is equal to
  $$
  \PAR{
    \begin{array}{rcl}
      (\Phi'''(u+v)-\Phi'''(u))\,v
      & \quad & \Phi''(u+v)-\Phi''(u) + \Phi'''(u+v)\,v \\
      \Phi''(u+v)-\Phi''(u) + \Phi'''(u+v)\,v
      & \, & 2\,\Phi''(u+v)+\Phi'''(u+v)\,v
    \end{array}},
  $$
  $\zeta_2$ is convex in $u$ but not necessarily in $v$.  Finally, let us
  show why \HYPOT{2'} implies \HYPOT{2}. Assume that \HYPOT{2'} holds.  The
  function $\Psi$ is non negative since $\Phi$ is convex, and it is convex in
  each variable $u$ and $v$ when $\Phi$ and $\Phi''$ are both convex.
  Moreover, it is bivariate convex when $\Phi(x)=x\log x$ or $\Phi(x)=x^2$, as
  observed in the pretty paper \cite{MR2002f:60109}.  Actually, this
  bi-convexity holds in much more cases since
  $$
  \GR^2\,\Psi(u,v):=
  \PAR{
    \begin{array}{rcl}
      \Phi''(u+v)-\Phi''(u)-\Phi'''(u)\,v & \quad & \Phi''(u+v)-\Phi''(u) \\
    \Phi''(u+v)-\Phi''(u)            & \quad &  \Phi''(u+v)
   \end{array}
 },
 $$
 for which $\mathrm{Tr}\,\GR^2\,\Psi(u,v) \geq 0$ when $\Phi$ and $\Phi''$ are 
 convex, and
 $$
 \mathrm{Det}\,\GR^2\,\Psi(u,v) =
 \Phi''(u)\,(\Phi''(u+v)-\Phi''(u))-\Phi''(u+v)\,\Phi'''(u)\,v.
 $$
 which is non negative since $\Phi$ and $\Phi''$ are convex and $\Phi'$ is
 concave. Alternatively, one can use Gershgorin-Hadamard Theorem (cf.
 \cite[Sect. 6.1]{MR91i:15001}) to see directly that $\GR^2\,\Psi$ is non 
 negative under \HYPOT{2'}.
\end{rem}

\section{Links with Boltzmann-Shannon entropy}
\label{se:shannon}

We assume here that $\cI=\dR_+$. Let $f$ be a probability density function on
$\dR^d$ with respect to the Lebesgue measure.  When $\Phi$ is convex, one can
define the Shannon $\Phi$-entropy of $f$ by:
\begin{equation}
  \HENT{}{\Phi}{f}
  :=-\ENT{}{\Phi}{f\,dx\,\vert\,dx}
  :=-\ENT{dx}{\Phi}{f}
  :=-\int_{\dR^d}{\whPhi(f)(x)\,dx},\label{eq:def-ent-sha-cont}
\end{equation}
where $\whPhi(u):=\Phi(u)-\Phi(1)\,u$. For any random vector $X$ in $\dR^d$
with density $f$ with respect to the Lebesgue measure, we denote
$\HENT{}{\Phi}{X}:=\HENT{}{\Phi}{f}$. This functional is invariant by
translations and hence does not depend on the mean of $f$.  We recover Shannon
entropy (or negentropy) $\HENTF{}{}$ when $\Phi(x)=x\log x$, cf. \cite[Chap.
10]{MR2002g:46132}.  Actually, $\SENTF{}{}:=-\HENTF{}{}$ is also known as
Shannon Information in Information Theory or Boltzmann entropy in Kinetic
Gases Theory, see for example \cite{MR50:9266} and \cite{MR1964483}.
Similarly, one can define the Shannon $\Phi$-entropy of a discrete probability
measure $p_1\de_{x_1}+\cdots+p_n\de_{x_n}$ by
\begin{equation}
  \HENT{}{\Phi}{p_1\de_{x_1}+\cdots+p_n\de_{x_n}}:=-\sum_{i=1}^n
  \whPhi(p_i).\label{eq:def-ent-sha-discr}
\end{equation}
Since $x\in\dR^n\mapsto\whPhi(x_1)+\cdots+\whPhi(x_n)$ is
convex, $\HENTF{}{\Phi}$ is a concave functional on the simplex
$$
\BRA{(p_1,\ldots,p_n) \in\dR_+\times\cdots\times\dR_+,\, p_1+\cdots+p_n=1}.
$$
At fixed $n$, it achieves its minimum $0$ for Dirac measures, which are the
extremal points of the simplex above, and its maximum $-n\,\whPhi(1/n)$ for the
uniform probability measure by convexity. Notice that the continuous version
is not always non-negative.

The important sub-additivity of Shannon entropy $\HENTF{}{}$ states that for
any random vector $(X_1,\ldots,X_n)$ with an absolute continuous law with
respect to the Lebesgue measure:
\begin{equation}
  \HENT{}{}{(X_1,\ldots,X_n)}\leq\HENT{}{}{X_1}+\cdots+\HENT{}{}{X_n},
  \label{eq:tenso-ent-sha}
\end{equation}
with equality if and only if $X_1,\ldots,X_n$ are independent. Such a
property relies on the non-negativity of Kullback-Leibler relative 
entropy and on
the basic additivity of the logarithm: $\log(ab)=\log a+\log b$:
$$
0\leq \ent{}{\cL((X_1,X_2))\,\vert\,\cL(X_1)\otimes\cL(X_2)}
= \HENT{}{}{(X_1,X_2)}-\HENT{}{}{X_1}-\HENT{}{}{X_2}.
$$
Notice that this sub-additivity is different from the one related to
Kullback-Leibler relative entropy \eqref{eq:ent-tenso} since Shannon
entropy is opposite in sign and based on the Lebesgue measure which is not a
probability measure. This fact was explained in \cite[Chap.
10]{MR2002g:46132}. Actually, such a sub-additivity property is not related to
the Lebesgue measure, as noticed in \cite{chafai-invlsg-2002}, and one can 
show that for any positive measure $\mu=\otimes_{i=1}^n \mu_i$ on a product 
space and any non-negative real valued integrable function $f$:
\begin{equation}
  \ent{\mu}{f} \geq
  \ent{\mu_1}{\int\!\!f\,d\mu_{\bs 1}}
  +\cdots+\ent{\mu_n}{\int\!\!f\,d\mu_{\bs n}},
  \label{eq:tenso-ent-sha-gen}
\end{equation}
where $\mu_{\bs i}:=\mu_1\times\cdots\times \mu_{i-1}\times
\mu_{i+1}\times\cdots\times\mu_n$, with equality if and only if $f$ is a
tensor product function. Beware that \eqref{eq:ent-tenso} and
\eqref{eq:tenso-ent-sha-gen} are opposite. We ignore if the sub-additivity
property \eqref{eq:tenso-ent-sha} of Shannon entropy can be generalised to any
convex $\Phi$. Actually, the functional $f\mapsto \HENT{}{\Phi}{f}$ is concave
and formally, the Fr\'{e}chet derivatives are given by:
$$
(D \HENTF{}{\Phi})(f)(h) = -\int_{\dR^d}\!\!\Phi'(f)\,h\,dx
\quad\text{and}\quad
(D^2 \HENTF{}{\Phi})(f)(h,h) = -\int_{\dR^d}\!\!\Phi''(f)\,h^2\,dx.
$$
Therefore, Shannon like $\Phi$-entropy $\HENTF{}{\Phi}{}$ achieves its
maximum under the linear constraint
$\moy{}{W(X)}=\int_{\dR^d}\!W(x)f(x)\,dx=c$ for probability densities
functions of the form:
$$
f_W := (\whPhi)'^{-1}\PAR{-\la -\be\,W},
$$
where $(\la,\be)\in\dR^2$ is chosen in such a way that the constraint is
fulfilled and that $f_W$ is a probability density function with respect to the
Lebesgue measure. Notice that $(\whPhi)'^{-1}$ is the derivative of the Young
conjugate of $\whPhi$. The Gaussian maximum of Shannon entropy $\HENTF{}{}$ at
fixed covariance appears as a particular case, for which $\Phi(x)=x\log x$,
$W(x)=\ABS{x}^2$, $(\whPhi)'^{-1}(y)=\exp(y-1)$, $\be=d/(2c)$ and
$\la=Z_{c,d}$ is the Gaussian normalising constant. More generally, we
recovered as a particular case the famous ``principle of maximum entropy''
which states that Boltzmann-Shannon entropy is maximised under linear
constraint by Boltzmann-Gibbs measures. It is tempting and quite natural to
ask if $\HENTF{}{\Phi}$ shares more common properties with Shannon entropy
$\HENTF{}{}$. As we have seen, the convex conjugate functional will play a
role. This question is partly answered in \cite{MR94b:94010},
\cite{MR87k:94007}, \cite{MR84e:62085} and \cite{MR84d:94009} and references
therein. Entropy like measure of information are still actively explored, and
one can find recent results in \cite{math.PR/0401368} and in Flemming
Tops{\o}e papers for example.
              
\bigskip

\textbf{Final words.} One may retain that the classical relative entropy
$\ent{}{\nu\,\vert\,\mu}$ where $\mu$ and $\nu$ are positive Borel measures
possesses a lot of properties coming from the very particular base function
$\Phi(x)=\Te(x):=x\log x$:
\begin{itemize}
\item $\Te$ is strictly convex, and thus $\ent{\mu}{f}=0$ iff $f$ is constant
  $\mu$-a.s.;
\item $\Te(0)=\Te(1)=0$ and thus $\bH(X)$ vanishes when $X$ is constant;
\item $1/\Te''$ is affine and hence concave and thus $f\mapsto \ent{\mu}{f}$
  is convex;
\item $\Te(ab)=b\,\Te(a)+a\,\Te(b)$ and thus $f\mapsto\ent{\mu}{f}$ is
  $1$-homogeneous;
\item The Young conjugate $\Te^*(u)=e^{u-1}$ is monotone and $\Te^{*'}=\Te^*$.
\end{itemize}
Recall that $\Te^*(u):=\int_0^u\Te'^{-1}(x)\,dx$. Some of these properties are
well imitated by $x\mapsto{}\ABS{x}^p$ with $p\in(1,2]$, which is exactly the
family of simple power convex functions between $x\mapsto\Te(x)$ and
$x\mapsto{}x^2$ and the latter is for some aspects the ``simplest'' one. Some
results involving $\entf{\mu}$ rely only on few properties of $\Te$ whereas
other ones rely on all of them, and this fact obviously puts some limits on
the possible generalisations.

\textbf{Acknowledgements.} This research was partially supported by MathFIT
(EPSRC \& LMS) postdoctoral grant GR/R2962/8/01 (HBKBU). The author would like
to thank Prof. Bogus{\l}aw Zegarli{\'n}ski for providing
\cite{bobkov-zegarlinski-02} and for interesting discussions, Prof. Franck
Barthe for providing a preliminary version of \cite{MR1863703}, and Prof.
Nicolas Privault and Prof. Feng-Yu Wang for communicating their preprints
\cite{MR2004c:60043} and \cite{wang-02} respectively.

\addcontentsline{toc}{section}{\refname}

\providecommand{\etalchar}[1]{$^{#1}$}
\providecommand{\bysame}{\leavevmode ---\ }
\providecommand{\og}{``}
\providecommand{\fg}{''}
\providecommand{\smfandname}{et}
\providecommand{\smfedsname}{\'eds.}
\providecommand{\smfedname}{\'ed.}
\providecommand{\smfmastersthesisname}{M\'emoire}
\providecommand{\smfphdthesisname}{Th\`ese}

\begin{center}
  \hrule
\end{center}

{\footnotesize\noindent%
  \textbf{Affiliations:} 
  \begin{itemize}
  \item UMR 181 de Physiopathologie et Toxicologie Exp\'erimentales INRA-ENVT,
    \'{E}cole Nationale V\'{e}t\'{e}-rinaire de Toulouse, 23 Chemin des
    Capelles, F-31076 Toulouse CEDEX 3, France.\\
    \textbf{E-mail:} \url{mailto:d(DOT)chafai(AT)envt(DOT)fr}\\
    \textbf{Web:} \url{http://biostat.envt.fr/~chafai/}
  \item Laboratoire de Statistique et Probabilit\'{e}s, UMR CNRS C05583,
    Institut de Math\'{e}matiques, Universit\'{e} Paul Sabatier, 118 route de
    Narbonne, F-31062 Toulouse CEDEX 4, France.\\
    \textbf{E-mail:} \url{mailto:chafai(AT)math(DOT)ups-tlse(DOT)fr}\\
    \textbf{Web:} \url{http://www.lsp.ups-tlse.fr/Chafai/}
  \item University of Oxford, Stochastic Analysis Group, Mathematical
    Institute, 24--29 St. Giles', Oxford OX1 3LB, United Kingdom.
\end{itemize}

\begin{center}
  \bigskip {--- \small Compiled \today{} ---}

  This paper is a reprint with minor corrections of the original paper published
  in Journal of Mathematics of Kyoto University, vol. 44 (2002), no. 2, 325-363.
\end{center}
}

\end{document}